\documentclass{amsart}
\usepackage{graphicx}
\usepackage{amssymb}
\usepackage{epstopdf}
\usepackage{nicefrac}
\usepackage{enumerate}
\usepackage{verbatim}
\usepackage{hyperref}
\usepackage{color}
\usepackage{pst-all}
\usepackage{tikz}

\newtheorem{thm}{THEOREM}[section]

\newtheorem{cor}[thm]{COROLLARY}

\newtheorem{defn}[thm]{DEFINITION}
\newtheorem{ex}[thm]{EXAMPLE}

\newtheorem{lemma}[thm]{LEMMA}
\newtheorem{prob}[thm]{PROBLEM}
\newtheorem{prop}[thm]{PROPOSITION}

\newtheorem{remark}[thm]{REMARK}

\newcommand{\ds}{\displaystyle}

\newcommand{\mC}{{\mathbb C}}

\newcommand{\mQ}{{\mathbb Q}}

\newcommand{\mZ}{{\mathbb Z}}
\newcommand{\mbZ}{{\mathbf Z}}

\newcommand{\cA}{{\mathcal A}}
\newcommand{\cB}{{\mathcal B}}

\newcommand{\cD}{{\mathcal D}}

\newcommand{\cF}{{\mathcal F}}
\newcommand{\cG}{{\mathcal G}}

\newcommand{\cK}{{\mathcal K}}
\newcommand{\cL}{{\mathcal L}}
\newcommand{\cM}{{\mathcal M}}
\newcommand{\cN}{{\mathcal N}}
\newcommand{\cO}{{\mathcal O}}
\newcommand{\cP}{{\mathcal P}}

\newcommand{\fX}{{\mathfrak{X}}}

\input xy
\xyoption {all}

\begin{document}

\title{Galois groups and Cantor actions}

\begin{abstract}
{ In this paper, we study the actions of profinite groups on Cantor sets which arise from representations of Galois groups of certain fields of rational functions. Such representations are associated to polynomials, and they are called profinite iterated monodromy groups. We are interested in a topological invariant of such actions called the asymptotic discriminant. 
In particular, we give a complete classification by whether the asymptotic discriminant is stable or wild in the case when the polynomial generating the representation is quadratic. We also study different ways in which a wild asymptotic discriminant can arise.

}
\end{abstract}

 \author{Olga Lukina}
 \email{olga.lukina@univie.ac.at}
\address{Faculty of Mathematics, University of Vienna, Oskar-Morgenstern-Platz 1, 1090 Vienna, Austria}

 \date{}


\thanks{Version date: September 21, 2018. Revision: May 15, 2020}

\date{}


\maketitle



\section{Introduction}\label{sec-intro}

In this paper, we consider actions of countable groups $G$, associated to representations of absolute Galois groups of number fields into automorphism groups of trees. Such representations are associated to polynomials of degree $d \geq 2$. Our main result is a complete classification of such actions by an invariant called the \emph{asymptotic discriminant} in the case $d = 2$. We also investigate the dynamical properties of such actions, and give criteria for determining the properties of the asymptotic discriminant which are applicable for these and more general actions.

Countable groups considered in this paper arise as dense subgroups of \emph{profinite iterated monodromy groups}. We now briefly recall necessary background on profinite iterated monodromy groups and the asymptotic discriminant in order to state our results.

\medskip
Let $K$ be a number field, that is, $K$ is a finite algebraic extension of the rational numbers $\mQ$. 
Let $f(x)$ be a polynomial of degree $d \geq 2$ with coefficients in the ring of integers of $K$. Let $t$ be a transcendental element, then $K(t)$ is the field of rational functions with coefficients in $K$.

Denote by $f^n(x)$ the $n$-th iterate of $f(x)$, and, for $n \geq 1$, consider the solutions of the equation $f^n(x) = t $ over $K(t)$. The polynomial $f^n(x)-t$ is separable and irreducible over $K(t)$ for all $n \geq 1$ \cite[Lemma 2.1]{AHM2005}. Therefore, it has $d^n$ distinct roots, and the Galois group $H_n$ of the extension $K_n$, obtained by adjoining to $K(t)$ the roots of $f^n(x)-t$, acts transitively on the roots. 

We represent the tower of preimages of $t$ under the iterations of $f$ as a tree $T$ as follows. Let $V_0$ be a singleton, called the \emph{root} of the tree $T$. For each $n \geq 1$, let the vertex set $V_n$ be the set of distinct roots of $f^n(x) = t$. Thus for $n \geq 1$ we have $|V_n| = d^n$. We join $ \beta \in V_{n+1}$ and $\alpha \in V_n$ by an edge if and only if $f(\beta) = \alpha$. For each $n \geq 1$, the Galois group $H_n$ acts transitively on the roots of $f^n(x) = t$ by field automorphisms, and so induces a permutation of vertices in $V_n$. The field extensions satisfy $K_{n} \subset K_{n+1}$, and so an automorphism of $V_{n+1}$ induces an automorphism of $V_n$, thus defining a group homomorphism $\lambda^{n+1}_n:H_{n+1} \to H_n$. Taking the inverse limit
   $${\rm Gal}_{\rm arith}(f) = \lim_{\longleftarrow}\{\lambda^{n+1}_n: H_{n+1} \to H_n\},$$
we obtain a profinite group called the \emph{arithmetic iterated monodromy group} of the polynomial $f(x)$. In other words, the group ${\rm Gal}_{\rm arith}(f)$ is the Galois group of the extension $\cK =  \bigcup_{n \geq 1} K_n $ obtained by adjoining to $K(t)$ the roots of $f^n(x) - t$ for $n \geq 1$. The action of the groups $H_n$ preserve the connectedness of the tree $T$, and so ${\rm Gal}_{\rm arith}(f)$ acts by permuting paths in $T$. We denote by $\cP_d$ the space of paths in a $d$-ary tree $T$. The space $\cP_d$ is topologically a Cantor set, see Section \ref{tree-model} for details.

Recall that $\overline{K}$ is a separable closure of $K$, and let $L = \overline{K} \cap \cK$ be the maximal constant field extension of $K$ in $\cK$, that is, $L$ contains all elements of $\cK$ algebraic over $K$. The Galois group ${\rm Gal}_{\rm geom}(f)$ of the extension $\cK/L(t)$, called the \emph{geometric iterated monodromy group}, is a normal subgroup of ${\rm Gal}_{\rm arith}(f)$. As we discuss in more detail in Section \ref{sec-IMG}, the group ${\rm Gal}_{\rm geom}(f)$ can be computed using methods from Geometric Group Theory. We outline these methods in Sections \ref{sec-IMG} and \ref{sec-wreathproduct}, and use them to prove one of our main theorems, Theorem \ref{thm-1}.
 
Both profinite groups ${\rm Gal}_{\rm arith}(f)$ and ${\rm Gal}_{\rm geom}(f)$ are inverse limits of finite groups, indexed by natural numbers, so by \cite[Proposition 4.1.3]{Wilson} they contain countably generated dense subgroups $G_{\rm arith}$ and $G_{\rm geom}$ respectively. By a slight abuse of notation, we denote by the same symbols the groups $G_{\rm arith}$ and $G_{\rm geom}$ with discrete topology. Thus associated to the actions of profinite monodromy groups ${\rm Gal}_{\rm arith}(f)$ and ${\rm Gal}_{\rm geom}(f)$, there are actions $(\cP_d, G_{\rm arith})$ and $(\cP_d, G_{\rm geom})$ of discrete groups on the Cantor set $\cP_d$. 
 
 \medskip
The \emph{asymptotic discriminant} is an invariant which classifies actions of discrete groups on Cantor sets. We now briefly introduce it and the property of actions which it detects. The reader should remember that the notion of the asymptotic discriminant is completely different to the notion of the `discriminant of a polynomial', which is the product of squares of differences of the polynomial roots. These two discriminants should not be confused. 

Let $X$ be a Cantor set, and let $G$ be a countably generated discrete group which acts on $X$ by homeomorphisms via the homomorphism $\Phi: G \to Homeo(X)$. We denote by $(X,G,\Phi)$ the action, and also use the notation $g \cdot x = \Phi(g)(x)$. We assume that the action is \emph{minimal}, that is, the orbit of any $x \in X$ under the action of $G$ is dense in $X$. We also assume that the action $(X,G,\Phi)$ is \emph{equicontinuous}, see Section \ref{tree-model} for a definition. Actions on trees in which we are interested in this paper are minimal and equicontinuous, see Section \ref{ex-treecylinder}. 

The concept of a \emph{quasi-analytic action} was first introduced by Haefliger  \cite{Haefliger1985} for the study of pseudogroups of local isometries on locally connected spaces. \'Alvarez L\'opez, Candel and Moreira Galicia \cite{ALC2009,ALM2016} reformulated Haefliger's definition to include totally disconnected spaces, and a local version of the definition as below first appeared in \cite{DHL2017}. 

\begin{defn}\label{defn-lqa-rev}
Let $G$ be a (countable or profinite) group. An action of $G$ on a metric space $X$ is \emph{\bf locally quasi-analytic (or LQA)} if there exists $\epsilon >0$ such that for any open set $U \subset X$ with ${\rm diam}(U)< \epsilon$ the following holds: for every $g_1,g_2 \in G$ and any open set $V \subset U$, if the restrictions $g_1|V = g_2|V$, then $g_1|U = g_2|U$. An action $(X,G,\Phi)$ is \emph{\bf quasi-analytic} if we can choose $U = X$ in this definition.
\end{defn}

If an action of $G$ on $X$ is LQA, then its elements have unique extensions from small open sets to open sets of diameter $\epsilon >0$. We use the LQA property to divide actions on Cantor sets into two large classes, \emph{stable} and \emph{wild} actions, as follows.

Since $(X,G,\Phi)$ is equicontinuous, then the closure of $\Phi: G \to Homeo(X)$ in the uniform topology on $Homeo(X)$ is a profinite group, identified with the Ellis (enveloping) group of the action \cite{Ellis1969,Auslander1988}, see also Section \ref{subsec-ellis}. The profinite group $\overline{\Phi(G)}$ acts on $X$ by homeomorphisms.

\begin{defn}\label{defn-stable-wild-rev}
Let $(X,G,\Phi)$ be a minimal equicontinuous group action, and let $\overline{\Phi(G)}$ be its Ellis group. Then $(X,G,\Phi)$ is \emph{\bf stable} if the action of the Ellis group $\overline{\Phi(G)}$ on $X$ is locally quasi-analytic (LQA), and the action $(X,G,\Phi)$ is \emph{\bf wild} otherwise.
\end{defn}

The image $\Phi(G)$ is a dense subgroup of $\overline{\Phi(G)}$, and the homomorphism $\Phi$ may have non-trivial kernel. In practice the subgroup $\Phi(G)$ is the easiest part of the Ellis group to work with. For example, in some cases one can show that the action of the countable group $\Phi(G)$ on $X$ is not LQA by finding for each open set $U \subset X$ elements which do not satisfy Definition \ref{defn-lqa-rev} on $U$. If the action of the dense subgroup $\Phi(G)$ is not LQA, then the action of its closure $\overline{\Phi(G)}$ is not LQA. It is not known at the moment if the converse implication holds. Indeed, it is conceivable that elements which do not have unique extensions are contained in $\overline{\Phi(G)} - \Phi(G)$, and so the action of $\Phi(G)$ is LQA while the action of $\overline{\Phi(G)}$ is not LQA.

The \emph{asymptotic discriminant} of a minimal equicontinuous group action $(X,G,\Phi)$ is a computable algebraic invariant which allows us to determine if an action of the Ellis group $\overline{\Phi(G)}$ is LQA, and so to determine if $(X,G,\Phi)$ is stable or wild. The work \cite{HL2017}, where this invariant was introduced, computes explicitly the asymptotic discriminant for the actions of a class of torsion free finite index subgroups of ${\rm SL}(n,\mZ)$. The rigorous definition of the asymptotic discriminant is quite technical, and so we postpone it until Section \ref{tree-model}, where we also discuss various examples and further applications. 

\begin{defn}
We say that an action $(X,G,\Phi)$ has \emph{\bf stable} (resp. \emph{\bf wild}) \emph{\bf asymptotic discriminant} if an action is stable (resp. wild) according to Definition \ref{defn-stable-wild-rev}.
\end{defn}

A stable asymptotic discriminant determines a closed subgroup of $\overline{\Phi(G)}$ which can be trivial, finite or an infinite profinite group. Thus we can speak about stable actions with trivial, finite, or infinite discriminant group, see Section \ref{tree-model} for details.

Wild actions may arise in different ways. For example, it was shown in the joint work with Hurder \cite{HL2018} that if a group contains a so-called \emph{non-Hausdorff element}, then its action is not LQA. 

\begin{defn}\label{defn-nonH-rev}
Let $G$ be a countable or profinite group acting on a topological space $X$. Then $g \in G$ is called a \emph{\bf non-Hausdorff element} if there is a point $x \in X$ such that $g \cdot x = x$, $g$ is not the identity map on any open neighborhood $W \owns x$, and for any such neighborhood $W$ there is an open set $O \subset W$ such that $g|O = id$. 
\end{defn}

The term \emph{non-Hausdorff element} is motivated by the topological properties of a germinal groupoid associated to the action, as we discuss in detail in Section \ref{subsec-nonHaus}. There are many examples of actions with non-Hausdorff elements, for example, three out of four generators of the Grigorchuk group in \cite[Section 1.6]{Nekr} are non-Hausdorff. Also, the process of the fragmentation of dihedral groups in \cite{Nekr2018} is done by adding non-Hausdorff elements to the group. As a consequence of our study of the asymptotic discriminant for actions of profinite iterated monodromy groups, we will obtain examples of non-LQA actions of discrete groups with no non-Hausdorff elements. We will also give a sufficient condition under which the closure $\overline{\Phi(G)}$ of the action contains a non-Hausdorff element.

\medskip
We now state our results. Let $f(x)$ be a polynomial of degree $d = 2$ with coefficients in the ring of integers of the field $K$, where $K$ is a finite extension of $\mQ$. The groups ${\rm Gal}_{\rm arith}(f)$ and ${\rm Gal}_{\rm geom}(f)$ are profinite groups, associated to the polynomial, which act on the space of paths $\cP_2$ of the binary tree $T$. The groups $G_{\rm arith}$ and $G_{\rm geom}$ are the respective countable dense subgroups of these profinite groups.

Denote by $c$ the critical point of $f(x)$, and by $P_c = \{f^m(c) \mid m \geq 1\}$ the orbit of the critical point. Polynomials in Theorems \ref{thm-1} and \ref{thm-3} are post-critically finite, which implies that $P_c$ is a finite set.

\begin{thm}\label{thm-1}
Let $K$ be a finite extension of $\mQ$, and let $f(x)$ be a quadratic polynomial with coefficients in the ring of integers of $K$, such that the orbit $P_c$ of the critical point $c$ of $f(x)$ is finite of length $r = \#P_c$. Consider the action of the profinite iterated geometric monodromy group ${\rm Gal}_{\rm geom}(f)$ on the space of paths $\cP_2$ of the tree $T$ of the solutions to $f^n(x) = t$, $n \geq 1$. Then the action has wild asymptotic discriminant, unless $r=1$ or the orbit of $c$ is strictly pre-periodic and $r = 2$.

More precisely, the following holds:
\begin{enumerate}
\item If the orbit of $c$ is strictly periodic, and $r = 1$, then the action of ${\rm Gal}_{\rm geom}(f)$ is LQA with trivial discriminant group.
\item If the orbit of $c$ is strictly periodic, and $r \geq 2$, then ${\rm Gal}_{\rm geom}(f)$ is conjugate in $Aut(T)$ to the closure of the action of a discrete group $\widetilde{G}_{r}$ on $r$ generators. The action of $\widetilde{G}_r$ is not LQA, and $\widetilde{G}_r$ does not contain non-Hausdorff elements. Thus the action of ${\rm Gal}_{\rm geom}(f)$ is not LQA and has wild asymptotic discriminant.
\item If the orbit of $c$ is strictly pre-periodic, and  $r = 2$, then the action of ${\rm Gal}_{\rm geom}(f)$ is LQA with finite discriminant group.
\item If the orbit of $c$ is strictly pre-periodic, and $r \geq 3$, then ${\rm Gal}_{\rm geom}(f)$ is conjugate in $Aut(T)$ to the closure of the action of a discrete group $\widetilde{H}_r$ on $r$ generators. The action of $\widetilde{H}_r$ is not LQA, and $\widetilde{H}_r$ contains non-Hausdorff elements. Thus the action of ${\rm Gal}_{\rm geom}(f)$ is not LQA and has wild asymptotic discriminant.

\end{enumerate}
\end{thm}
The groups $\widetilde{G}_r$, $r \geq 2$, and $\widetilde{H}_r$, $r \geq 3$, in Theorem \ref{thm-1} are described in detail in Section \ref{sec-geom}.

Polynomials in case $(1)$ of Theorem \ref{thm-1} are conjugate to $x^2$, and polynomials in case $(3)$ are conjugate to the Chebyshev polynomial $x^2-2$. Thus Theorem \ref{thm-1} states that if $f$ is quadratic, then the action of ${\rm Gal}_{\rm geom}(f)$ has wild asymptotic discriminant, unless $f$ is conjugate to a Chebyshev polynomial or a powering map. 

If the orbit $P_c$ of the critical point $c$ of $f(x)$ is infinite, then by \cite[Section 1.10]{Pink2013} ${\rm Gal}_{\rm geom}(f) \cong Aut(T)$, and so  by \cite{Lukina2018} the action of ${\rm Gal}_{\rm geom}(f)$ on the tree of solutions to $f^n(x) = t$, $n \geq 1$, is not LQA and has wild asymptotic discriminant. This together with the results of Theorem \ref{thm-1} provides a complete classification by the asymptotic discriminant of the actions of profinite geometric iterated monodromy groups associated to quadratic polynomials defined over finite extensions of $\mQ$. 

The case $(2)$ of Theorem \ref{thm-1} is the case when the critical point $c$ of a quadratic polynomial $f(x)$ has a finite strictly periodic orbit of cardinality at least $2$. The groups $\widetilde{G}_r$, $r \geq 2$, in $(2)$ give rise to a family of actions of discrete groups which are not LQA and where the groups do not contain non-Hausdorff elements. These are the first examples of this kind, known to the author. The closures of the actions of $\widetilde{G}_r$ may or may not contain non-Hausdorff elements, see Sections \ref{subsec-conjugacy} and \ref{sec-geom} for more discussion. We summarize this discussion as a corollary of Theorem \ref{thm-1}.

\begin{cor}\label{cor-1}
There exist actions of discrete finitely generated groups on a Cantor set $X$ which are not LQA and such that the groups do not contain non-Hausdorff elements. That is, non-LQA minimal Cantor actions of discrete groups with and without non-Hausdorff elements form two distinct non-empty classes.
\end{cor}

The following problem remains open.

\begin{prob}\label{prob-closurenonHausdorff}
{\rm
Let $\Phi:G \to Homeo(X)$ be an action of a countable group $G$ on a Cantor set $X$. Suppose that the action is not LQA. Show that the closure $\overline{\Phi(G)}$ contains a non-Hausdorff element, or find a counterexample.
}
\end{prob}

Another subtle point of the statement $(2)$ and $(4)$ in Theorem \ref{thm-1} is that although ${\rm Gal}_{\rm geom}(f)$ is conjugate in $Aut(T)$ to the closure of the action of $\widetilde{G}_r$ and, using the conjugacy and Proposition \ref{prop-conjugacy}, we can determine that the action of ${\rm Gal}_{\rm geom}(f)$ is not LQA, we cannot say much about the action of a dense subgroup $G_{\rm geom}$ of ${\rm Gal}_{\rm geom}(f)$. Indeed, $\widetilde{G}_r$ need not be mapped onto $G_{\rm geom}$ under the conjugacy. Nekrashevych \cite{Nekr2007} gives examples of actions of discrete non-isomorphic groups on $3$ generators whose closures in $Aut(T)$ are conjugate. Since the asymptotic discriminant is an invariant of the closures of the actions, we conclude that the action is wild without knowing if the action of $G_{\rm geom}$ is LQA or not LQA.

Polynomials described by statement $(3)$ in Theorem \ref{thm-1} are examples of quadratic polynomials conjugate to Chebyshev polynomials. 
Chebyshev polynomials $T_d$ of degree $d \geq 2$ are described in Section \ref{sec-cheb}. These polynomials are well-studied, with discrete iterated monodromy groups ${\rm IMG}(T_d)$ computed for all degrees $d \geq 2$. As explained in Section \ref{sec-IMG}, the closure of the group ${\rm IMG}(T_d)$ is isomorphic to the profinite geometric iterated monodromy group of $T_d$. We have the following theorem. 

\begin{thm}\label{thm-cheb}
Let $T_d$ be the Chebyshev polynomial of degree $d \geq 2$ over $\mC$. Then the action of ${\rm IMG}(T_d)$ is stable with discriminant group $\mZ/2\mZ$, the finite group of order $2$.
\end{thm}

Let us now consider the action of the arithmetic iterated monodromy group ${\rm Gal}_{\rm arith}(f)$ for quadratic polynomials over $K$, where $K$ is a finite extension of $\mQ$. Recall that ${\rm Gal}_{\rm geom}(f)$ is a normal subgroup of ${\rm Gal}_{\rm arith}(f)$, and so ${\rm Gal}_{\rm arith}(f)$ is a subgroup of the normalizer $N$ of ${\rm Gal}_{\rm geom}(f)$ in $Aut(T)$. Properties of the normalizer of ${\rm Gal}_{\rm geom}(f)$, and how ${\rm Gal}_{\rm arith}(f)$ sits in the normalizer were studied by Pink \cite{Pink2013}. Direct computations based on the results of \cite{Pink2013} in $(1)$ and $(3)$, and Lemma \ref{lemma-inclusionwild} of this paper  in $(2)$ and $(4)$ give the following theorem.

\begin{thm}\label{thm-3}
Let $K$ be a finite extension of $\mQ$, and let $f(x)$ be a quadratic polynomial with coefficients in the ring of integers of $K$, such that the orbit $P_c$ of the critical point $c$ of $f(x)$ is finite of length $r = \#P_c$.  Consider the action of the profinite arithmetic iterated monodromy group ${\rm Gal}_{\rm arith}(f)$ on the path space $\cP_2$ of the tree $T$ of the solutions to $f^n(x) = t$, $n \geq 1$. Then the action has wild asymptotic discriminant, unless $r=1$ or the orbit of $c$ is strictly pre-periodic and $r = 2$. 

More precisely, the following holds:
\begin{enumerate}
\item If the orbit of $c$ is strictly periodic, and $r = 1$, then the action of ${\rm Gal}_{\rm arith}(f)$ is stable with infinite discriminant group.
\item If the orbit of $c$ is strictly periodic, and $r \geq 2$, then the action of ${\rm Gal}_{\rm arith}(f)$ is not LQA  and so the asymptotic discriminant of the action is wild. 
\item If the orbit of $c$ is strictly pre-periodic, and  $r = 2$, then the action of ${\rm Gal}_{\rm arith}(f)$  is stable with infinite discriminant group.
\item If the orbit of $c$ is strictly pre-periodic, and $r \geq 3$, then the action of ${\rm Gal}_{\rm arith}(f)$  is not LQA  and so  the asymptotic discriminant of the action is wild. 
\end{enumerate}
\end{thm}

If the critical point $c$ has infinite orbit, the action of ${\rm Gal}_{\rm geom}(f)$ is not LQA, and so Lemma \ref{lemma-inclusionwild} implies that the action of ${\rm Gal}_{\rm arith}(f)$ is not LQA and has wild asymptotic discriminant. This together with Theorem \ref{thm-3} completes the classification by the asymptotic discriminant of actions of profinite arithmetic iterated monodromy groups for quadratic polynomials over finite extensions of $\mQ$. 

A similar classification as in Theorems \ref{thm-1} and \ref{thm-3} for polynomials of degree $d \geq 3$ is currently out of our reach, mostly due to the absence of such comprehensive study of the geometric and arithmetic iterated monodromy groups in this case as was done by Pink for quadratic polynomials in \cite{Pink2013}. In some cases the question if a profinite iterated monodromy group, or its specialization, called an arboreal representation (see Section \ref{sec-IMG} for details), has wild asymptotic discriminant can be answered using the following algebraic criterion. The criterion is a sufficient condition for a profinite group acting on the path space $\cP$ of a spherically homogeneous rooted tree $T$ (not necessarily $d$-ary, see Section \ref{ex-treecylinder}) to contain a non-Hausdorff element. Actions with non-Hausdorff elements have wild asymptotic discriminant by the results of \cite{HL2018}.

\begin{thm}\label{thm-4}
Let $T$ be a spherically homogeneous rooted tree, and $\cP$ be the space of paths in  $T$. Let ${\ds H_\infty = \lim_{\longleftarrow}\{H_{n+1} \to H_n\}}$ be a profinite group, acting on $\cP$, so that for each $ n\geq 1$ the action of finite groups $H_n$ on the vertex sets $V_n$ of $T$ is transitive.  Let $\{L_n\}_{n \geq 1}$ be a collection of non-trivial finite groups such that for each $n \geq 1$ the group $H_n$ contains the wreath product $\cL_n = L_n \rtimes L_{n-1} \rtimes \cdots \rtimes L_1$. Then $H_\infty$ contains a non-Hausdorff element.
\end{thm}

Directly applying Theorem \ref{thm-4}, we conclude that the actions of the geometric and arithmetic iterated monodromy groups for the polynomial $f(z) = -2z^3 + 3z^2$ studied in \cite{BFHJY} have wild asymptotic discriminants, and a similar statement holds for a large class of dynamical Belyi maps with exactly three ramification points, studied in \cite{BEK2018}, see Example \ref{ex-11} for details.
A criterion under which the hypothesis of Theorem \ref{thm-4} holds for the arithmetic iterated monodromy group of a rational function $f(x)$ of degree $d \geq 2$ was proved in \cite[Theorem 3.1]{JKMT2015}. One of the conditions is that critical points of $f$ have non-intersecting orbits. For quadratic rational functions with two non-intersecting post-critical orbits it was shown in \cite[Theorem 4.8.1]{Pink2013-2}, that both the arithmetic and the geometric iterated monodromy groups are equal to $Aut(T)$ and so their actions are not LQA.

\medskip
We finish the introduction with a discussion of the motivation for the study of the asymptotic discriminant for actions arising from representations of Galois groups, and of the possible directions of future work.

One of the problems which motivated the study of arboreal representations in arithmetic dynamics is the problem of density of prime divisors in non-linear relations $a_n = f(a_{n-1})$. More precisely, let $f(x)$ be a polynomial of degree $d \geq 2$ with coefficients in a ring of integers of a field $K$, and let $a_0$ be a point in $K$. Let $Y_\infty$ be a representation of the absolute Galois group of $K$ into the group of automorphisms $Aut(T)$ of a $d$-ary tree $T$, as in Remark \ref{remark-specialization}. Recall that $\cP_d$ denotes the space of paths of the tree $T$. Consider the orbit $\cO(a_0) = \{f^n(a_0)\}_{n \geq 0}$.
What is the natural density of prime divisors of the points in $\cO(a_0)$? Odoni \cite{Odoni1985} showed that an upper estimate for such density can be obtained by counting the proportion $\cF(Y_\infty)$ of elements in $Y_\infty$ which fix at least one point in the space of paths $\cP_d$. The proportion $\cF(Y_\infty)$ is computed as the limit of the proportions of the elements with fixed points in the Galois groups $Y_n$ of finite extensions of $K$, obtained by adjoining the roots of the $n$-th iterate $f^n(x) = a_0$ to $K$. 
In particular, Odoni \cite{Odoni1985} showed that if $Y_\infty$ is isomorphic to the infinite wreath product $[S_d]^\infty$, where $S_d$ is a permutation group on $d$ elements, then $\cF(Y_\infty) = 0$. Jones \cite{Jones2008,Jones2013,Jones2015} developed a method of computing $\cF(Y_\infty)$ using theory of stochastic processes, and, in a series of papers, obtained the values of $\cF(Y_\infty)$ for various classes of arboreal representations. 

The discriminant group $\cD_x$ of an action $(X,G,\Phi)$ with profinite enveloping group $\overline{\Phi(G)}$ counts elements of $\overline{\Phi(G)}$ which fix a given point $x$. For every other point $y \in X$, the cardinality of $\cD_y$ is equal to that of $\cD_x$. Although the relationship between the cardinality of the asymptotic discriminant and the proportion $\cF(Y_\infty)$ is not direct, since they consider essentially the same objects, it is natural to ask if there is a relation between them. If an action has wild discriminant, then the number of elements in the finite Galois groups $Y_n$, whose cardinality is the denominator in the sequence computing the proportion $\cF(Y_\infty)$, grows extremely fast with $n$, so the following question is natural.

\begin{prob}
Let $Y_\infty$ be an image of an arboreal representation of a Galois group of a field, that is, a profinite group acting on the space of paths of a $d$-ary tree $T$. Suppose the action of $Y_\infty$ has wild asymptotic discriminant. Prove that the proportion $\cF(Y_\infty)$ of elements with fixed points in $Y_\infty$ is zero, or give a counterexample. 
\end{prob}  

One can also ask about the converse of this statement in the case when $Y_\infty$ has elements which fix points other than the identity, and when $\cF(Y_\infty) = 0$. If the only element in $Y_\infty$ which fixes points is the identity, then the action of $Y_\infty$ is stable.

Another motivation to consider the actions associated to representations of Galois groups into the automorphism groups of $d$-ary trees comes from the topological point of view. Actions of profinite iterated monodromy groups and of arboreal representations present a large class of examples, and studying them one can gain insights in the properties of general group actions on Cantor sets. For example, in this paper we studied the different ways in which a wild action can arise. 

As a consequence of Theorem \ref{thm-1}, in Corollary \ref{cor-1} we obtained that non-LQA minimal equicontinuous actions of discrete countable groups on Cantor sets may give rise to Hausdorff or non-Hausdorff \'etale groupoids $\cG(X,G,\Phi)$. It is interesting to compare these results to those obtained in related settings. For example, Hughes \cite{Hughes2012} considers germinal groupoids of \emph{all} local isometries on ultrametric spaces, which include Cantor sets. An interesting property of these groupoids is that if there is a point $x \in X$ and a local isometry $\ell$ of $X$ such that $\ell(x) = x$, then there exists a local isometry $\tilde{\ell}$ with $\tilde{\ell}(x) = x$ which is in our terminology a non-Hausdorff element. A consequence of this, in particular, is that ultrametric spaces with Hausdorff groupoids only admit local isometries where fixed points have fixed clopen neighborhoods. Hughes's setting is very different to ours, since in our setting every local map of the space $X$ must arise as a restriction of an action of an element $g \in G$. Still, it would be interesting to find out to what extent the properties of groupoids in these two different settings mirror each other.

The rest of the paper is organized as follows. In Section \ref{sec-IMG} we give more details about profinite and discrete iterated monodromy groups. In Section \ref{tree-model} we recall the method of group chains, and use it to introduce the asymptotic discriminant. In Section \ref{sec-wreathproduct} we recall the necessary background on wreath products and actions of self-similar groups as in \cite{Nekr}, and study the properties of non-Hausdorff elements in contracting groups. In Section \ref{sec-cheb} we prove Theorem \ref{thm-cheb} by an explicit computation of the asymptotic discriminant. The proof of Theorem \ref{thm-1} in Section \ref{sec-geom} relies on the explicit expressions for generators of $\widetilde{G}_r$ and $\widetilde{H}_r$, obtained by Pink \cite{Pink2013}. Each generator is associated to a point in the post-critical orbit and is defined recursively, so that how it acts on the tree depends on the order of the point in the critical orbit. For strictly periodic post-critical orbits, generators which are further away from the critical point act trivially on certain clopen subsets in $\cP_d$, and we use that to show that the action of ${\rm Gal}_{\rm geom}(f)$ on $\cP_d$ is not LQA. If the orbit is strictly pre-periodic, generators corresponding to the periodic part of the orbit are non-Hausdorff elements, and the action of ${\rm Gal}_{\rm geom}(f)$ is again not LQA.
Theorem \ref{thm-3} is proved in Section \ref{sec-arith}. By Lemma \ref{lemma-inclusionwild} if the action of ${\rm Gal}_{\rm geom}(f)$ is not LQA, then the action of ${\rm Gal}_{\rm arith}(f)$ is not LQA, so the main part of the proof of Theorem \ref{thm-3} is devoted to the cases where the polynomial $f$ is conjugate to the powering map or to a quadratic Chebyshev polynomial. In this cases the proof proceeds by explicitly computing the asymptotic discriminant of the action, using the information about the inclusion ${\rm Gal}_{\rm geom}(f) \subset {\rm Gal}_{\rm arith}(f)$ from \cite{Pink2013}.
The proof of Theorem \ref{thm-4}, which gives a condition under which a profinite group contains a non-Hausdorff element, is constructive and it given in Section \ref{sec-subgroups}.

\section{Iterated monodromy groups}\label{sec-IMG}

We recall some background about iterated monodromy groups.  In our description of profinite iterated arithmetic and geometric monodromy groups we follow \cite{Jones2013}, also see \cite{Pink2013}. 

Let $K$ be a number field, that is, $K$ is a finite algebraic extension of the rational numbers $\mQ$. 
Let $f(x)$ be a polynomial of degree $d \geq 2$ with coefficients in the ring of integers of $K$. Let $t$ be a transcendental element, then $K(t)$ is the field of rational functions with coefficients in $K$. 

We define the profinite iterated arithmetic and geometric monodromy groups as in the Introduction. That is, $f^n(x)$ denotes the $n$-th iterate of $f(x)$, and we consider the solutions of the equation $f^n(x) = t $ over $K(t)$. The polynomial $f^n(x)-t$ is separable and irreducible over $K(t)$ for all $n \geq 1$ \cite[Lemma 2.1]{AHM2005}. Therefore, it has $d^n$ distinct roots, and the Galois group $H_n$ of the extension $K_n$ obtained by adjoining to $K(t)$ the roots of $f^n(x)-t$ acts transitively on the roots. 

The tree $T$ has the vertex set $V = \bigsqcup_{n \geq 0} V_n$, where $V_0$ is identified with $t$, and $V_n$ with the sets of solutions of $f^n(x) = t$. We join $ \beta \in V_{n+1}$ and $\alpha \in V_n$ by an edge if and only if $f(\beta) = \alpha$. For each $n \geq 1$, the Galois group $H_n$ acts transitively on the roots of $f^n(x) = t$ by field automorphisms, and so induces a permutation of vertices in $V_n$. Since the field extensions satisfy $K_{n} \subset K_{n+1}$, we have a group homomorphism $\lambda^{n+1}_n:H_{n+1} \to H_n$. Taking the inverse limit
   $${\rm Gal}_{\rm arith}(f) = \lim_{\longleftarrow}\{\lambda^{n+1}_n: H_{n+1} \to H_n\},$$
we obtain a profinite group called the \emph{arithmetic iterated monodromy group} of the polynomial $f(x)$. For $n \geq 1$, the action of $H_n$ on $T$ preserves the connectedness of paths in $T$, and so ${\rm Gal}_{\rm arith}(f)$ is identified with a subgroup of the automorphism group $Aut(T)$ of the tree $T$. 

\begin{remark}\label{remark-specialization}
{\rm
Given a polynomial $f(x)$ over $K$, one can also consider the extensions $K(f^{-n}(\alpha))/K$ with Galois groups $Y_n$ for some $\alpha \in K$. If all iterates $f^n(x) - \alpha$ are separable and irreducible, by a similar procedure as above one can construct an  \emph{arboreal representation} ${\ds Y_\infty = \lim_{\longleftarrow}\{Y_{n+1} \to Y_n\}}$
  of the absolute Galois group of $\overline{K}/K$, where $\overline{K}$ is a separable closure of $K$, into the automorphism group $Aut(T_\alpha)$ of a tree $T_\alpha$. Since $f^n(x) - \alpha$ is irreducible for $n \geq 1$, the group $Y_n$ acts transitively on the vertex set $V_{\alpha,n}$ in $T_\alpha$. Since $f^n(x) - \alpha$ is separable, for each $n \geq 1$, $|V_{\alpha, n}| = |V_n| = d^n$, and $T_\alpha$ is a $d$-ary tree. It follows that $Aut(T)$ and $Aut(T_\alpha)$ are isomorphic, and  one may think of the groups $Y_n$ as obtained via the \emph{specialization} $t = \alpha$. As explained in \cite{Odoni1985}, Galois groups of polynomials do not increase under such specializations, and certain groups are preserved. So one can think of ${\rm Gal}_{\rm arith}(f)$ as $Y_\infty$ for a generic choice of $\alpha$, in a loose sense \cite{Odoni1985,Jones2013}.
 }
\end{remark}

Recall that $\overline{K}$ is a separable closure of $K$, and let $L = \overline{K} \cap \cK$ be the maximal constant field extension of $K$ in $\cK$, that is, $L$ contains all elements of $\cK$ algebraic over $K$. The Galois group ${\rm Gal}_{\rm geom}(f)$ of the extension $\cK/L(t)$ is a normal subgroup of ${\rm Gal}_{\rm arith}(f)$, and there is an exact sequence \cite{Pink2013,Jones2013}
 \begin{align}\label{exact-1} \xymatrix{ 1 \ar[r] & {\rm Gal}_{\rm geom}(f) \ar[r] &  {\rm Gal}_{\rm arith}(f) \ar[r] & {\rm Gal}(L/K) \ar[r] & 1}. \end{align}
The profinite group ${\rm Gal}_{\rm geom}(f)$ is called the \emph{geometric monodromy group}, associated to the polynomial $f(x)$. The geometric monodromy group ${\rm Gal}_{\rm geom}(f)$ does not change under extensions of $L$, so one can calculate ${\rm Gal}_{\rm geom}(f)$ over $\mC(t)$. 

Let $\mathbb{P}^1(\mathbb{C})$ be the projective line over $\mathbb{C}$ (the Riemann sphere), and extend the map $f: \mathbb{C} \to \mathbb{C}$ to $\mathbb{P}^1(\mathbb{C})$ by setting $f(\infty) = \infty$. Then $\infty$ is a critical point of $f(x)$. Let $C$ be the set of all critical points of $f(x)$, and let $P_{C} = \bigcup_{n \geq 1} f^n(C)$ be the set of the forward orbits of the points in $C$, called the \emph{post-critical set}. Suppose $P_{C}$ is finite, then the polynomial $f(x)$ is called \emph{post-critically finite}. 

If the polynomial $f(x)$ is post-critically finite, then it defines a partial $d$-to-$1$ covering $f: \cM_1 \to \cM$, where $\cM = \mathbb{P}^1(\mathbb{C}) \backslash P_C$ and $\cM_1 = f^{-1}(\cM)$ are punctured spheres. An element $s \in \cM$ has $d$ preimages under $f$, and $d^n$ preimages under the $n$-th iterate $f^n$. Denote by $\widetilde{V}_n = f^{-n}(s)$, for $n \geq 1$. In a manner similar to the one used to define the profinite iterated monodromy groups, one constructs a rooted $d$-ary tree $\widetilde{T}$ of preimages of $s$, with vertex sets $\bigsqcup_{n\geq 0} \widetilde{V}_n$. The fundamental group $\pi_1(\cM, s)$ acts on the vertex sets  of $\widetilde{T}$ via path-lifting. Let ${\rm Ker}$ be the subgroup of $\pi_1(\cM,s)$ consisting of elements which act trivially on \emph{every} vertex set $\widetilde{V}_n$, $n \geq 1$. The quotient group ${\rm IMG}(f) = \pi_1(\cM,s)/{\rm Ker}$, called the \emph{discrete iterated monodromy group} associated to the partial self-covering $f: \cM_1 \to \cM$, acts on the space of paths of the tree $\widetilde{T}$. If $f(x)$ is post-critically finite, by \cite[Proposition 6.4.2]{Nekr} attributed by Nekrashevych to R. Pink, ${\rm Gal}_{\rm geom}(f)$ over $\mC(t)$ is isomorphic to the closure of the action of ${\rm IMG}(f)$ in $Aut(\widetilde{T}) \cong Aut(T)$, where $T$ is the $d$-ary tree equipped with the action of ${\rm Gal}_{\rm geom}(f)$.

The actions of discrete iterated monodromy groups ${\rm IMG}(f)$ are well-studied, with many results known and many techniques developed, see, for example, Nekrashevych \cite{Nekr}. We outline some of these methods, needed for the proof of Theorem \ref{thm-1}, in Section \ref{sec-wreathproduct}.

\section{The asymptotic discriminant of an equicontinuous Cantor action}\label{tree-model}

In this section, we recall the necessary background on equicontinuous Cantor actions and the asymptotic discriminant. Main references for this section are works \cite{FO2002,DHL2016,DHL2017,HL2017}. Although the standing assumption in \cite{DHL2016,DHL2017,HL2017} was that $G$ is a finitely generated group, the reason for that was not any restrictions imposed by the proofs or by the properties of the objects considered. The motivation in those papers was to study and classify the dynamics of weak solenoids, and for a group to act on a Cantor fibre of a weak solenoid it must be realizable as a homomorphic image of a fundamental group of a closed manifold. Thus we assumed that the groups under consideration were finitely generated. However, the notion of the Ellis (enveloping) group does not require finite generation, and finite generation was not used in any of the proofs in \cite{DHL2016,DHL2017,HL2017}, so one easily checks that the results apply for the countably generated groups as well. 

\subsection{Equicontinuous actions on Cantor sets} 

Let $X$ be a Cantor set, that is, a compact totally disconnected metric space without isolated points. Recall that a space $X$ without isolated points is called \emph{perfect} \cite[Section 30]{Willard}.  

Let $D$ be a metric on $X$, and suppose $\Phi: G \to {Homeo}(X)$ defines an action of a countably generated discrete group $G$ on $X$. The action $(X,G,\Phi)$ is \emph{equicontinuous}, if for any $\epsilon >0$ there exists $\delta>0$ such that for any $g \in G$ and any $x,y \in X$ such that $D(x,y) < \delta$ we have $D(\Phi(g)(x),\Phi(g)(y))< \epsilon$.

All actions considered in this paper are minimal equicontinuous actions on Cantor sets. If an acting group is discrete, we call such an action a \emph{group Cantor action}.

\subsubsection{Boundary of a spherically homogeneous tree is a Cantor set}\label{ex-treecylinder}
For readers with non-topological background we explain why the set of infinite paths in a tree $T$ is a Cantor set, and why the path-preserving level-transitive action of $G$ is minimal and equicontinuous. 

A \emph{spherical index} is a sequence ${ \ell} = (\ell_1,\ell_2,\ldots)$ of positive integers, where $\ell_n \geq 2$ for $n \geq 1$. Let $T$ be a tree, defined by ${ \ell}$. That is, the set of vertices is $V = \bigsqcup_{n \geq 0} V_n$, where $V_0$ is a singleton, and  for $n \geq 1$  $V_n$ contains $\ell_1 \ell_2 \cdots \ell_n$ vertices. Since $|V_0|=1$, the tree $T$ is \emph{rooted}. Every vertex in $V_n$ is connected by edges to precisely $\ell_{n+1}$ vertices in $V_{n+1}$, and every vertex in $V_{n+1}$ is connected by an edge to precisely one vertex in $V_n$. A tree with this property is called a \emph{spherically homogeneous} tree. A path in $T$ is an infinite sequence $(v_n)_{n \geq 0} = (v_0,v_1,v_2,\cdots)$ such that $v_{n}$ and $v_{n+1}$ are connected by an edge, for $n \geq 0$. Denote by $\cP$ the set of all such sequences. If $\ell_n = d$ for all $n \geq 1$, where $d \geq 2$, then we call $T$ a \emph{$d$-ary rooted tree}, and denote the space of paths in $T$ by $\cP_d$. 

We are now going to define a topology on $\cP$ so that it is a Cantor set. 

For each $n \geq 0$, let $R_n  = \{0,1,\cdots, \ell_n-1\}$. If $T$ is $d$-ary, then for all $n \geq 1$ $R_1 = R_n$, and the set $R_1$ with $d$ elements is called the \emph{alphabet}. For $n \geq 1$, let ${\rm pr}_{n-1}:R_1 \times R_2 \times \cdots \times R_n \to R_1 \times R_2 \times \cdots \times R_{n-1}$ be the projection onto the product of the first $n-1$ sets. Define the bijections 
  \begin{align}\label{eq-labeling} b_n : R_1 \times R_2 \times \cdots \times R_n \to V_n \end{align}
 in such a way that 
  $${\rm pr}_{n-1}(b_n^{-1}(v)) = {\rm pr}_{n-1}(b_n^{-1}(w))$$ if and only if $v$ and $w$ are connected by edges to the same vertex in $V_{n-1}$.
Given a vertex $v \in V_n$, the preimage ${b}_n^{-1}(v) = (t_1,t_2,\cdots,t_n) $ is an $n$-tuple of integers, which we write as a word, that is, ${b}_n^{-1}(v) = t_1t_2\cdots t_n$. The collection of mappings \eqref{eq-labeling}, for $n \geq 1$, assigns a label to every vertex in the set $\bigcup_{n \geq 1} V_n$. We omit the label on the root as it is not used.

It follows from the definition of the maps \eqref{eq-labeling} that every word $t_1 t_2 \cdots t_n$ defines a finite path $(v_k)_{0 \leq k \leq n}$, where $v_0$ is the unique vertex in $V_0$, and $v_k = {b}_k(t_1 t_2 \cdots t_k)$ for $1 \leq k \leq n$. Then every infinite word $t_1t_2 \cdots t_n \cdots$, where $0 \leq t_n \leq \ell_n-1$, defines an infinite path in $T$, and so there is a bijection
  \begin{align}\label{eq-binfty} {b}_\infty: \prod_{n \geq 1} R_n  \to \cP \end{align}
such that ${b}_\infty|_{R_1 \times \cdots \times R_n} = {b}_n$. For $n \geq 0$, give $R_n$ discrete topology, then the product $\prod_{n \geq 1} R_n$ is compact by the Tychonoff theorem \cite{Willard}. Points are the only connected components in $\prod_{n \geq 1} R_n$, so $\prod_{n \geq 1} R_n$ is totally disconnected. 

Open sets in the product topology on $\prod_{n \geq 1} R_n$ have the form $\prod_{n \geq 0} U_n$, where $U_n \subseteq R_n$, and for all but a finite number of $n$ we have $U_n = R_n$. For example, given a word $t_1t_2 \cdots t_k$, let $U_n = \{t_n\}$ for $1 \leq n \leq k$, and $U_n = R_n$ otherwise. Then $U = \prod_{n \geq 1} U_n$ is the set of all infinite sequences in $\prod_{n \geq 1} R_n$ which start with a finite word $t_1 t_2 \cdots t_k$. In $\cP$ this set corresponds to all paths which contain the vertex $v_k = {b}_k(t_1t_2 \cdots t_k)$. We denote such a set $U$ by $U_k(t_1t_2 \cdots t_k)$.  

Note that for any open set $U=\prod_{n \geq 1} U_n$ its complement in $\prod_{n \geq 1} R_n$ is also open, and so $U$ is closed. A set which is open and closed is called a \emph{clopen set}. 

Let $t = t_1 t_2 \cdots $ be an infinite sequence, and consider a descending sequence of open neighborhoods $U_n(t_1 t_2 \cdots t_n)$ for $n \geq 1$. Since $|R_n| \geq 2$ for all $n \geq 1$, every $U_n(t_1 t_2 \cdots t_n)$ is infinite, and it follows that $\prod_{n \geq 1} R_n$ is a perfect set. 

We have shown that $\prod_{n \geq 1}R_n$, and so $\cP$, is a Cantor set. From now on we identify $\prod_{n \geq 1}R_n$ and $\cP$, and think of elements in $\cP$ as infinite sequences $t_1 t_2 \cdots$, where $0 \leq t_n \leq \ell_n - 1$ for $n \geq 1$. For $v_n = b_n(t_1t_2 \cdots t_n)$, we suppress $b_n$ in the notation, and just write $v_n = t_1 \cdots t_n$.

Let $G$ be a countably generated discrete group, and let $G$ act on the tree $T$ by permuting vertices in each $V_n$, $n \geq 0$, in such a way that the connectedness of paths in $T$ is preserved, and the action is transitive on each $V_n$. Since permutations are bijective, the action of each $g \in G$ induces a bijective map $\Phi(g):\cP \to \cP$. For each $n \geq 1$, the image of an open set $U_n(t_1 t_2 \cdots t_n)$ under $\Phi(g)$ is an open set $U_n(g \cdot (t_1 t_2 \cdots t_n))$, so $\Phi(g)$ is a homeomorphism. Thus $G$ acts on $\cP$ by homeomorphisms.

Let $w = w_1 w_2 \cdots \in \cP$ be an infinite sequence. Since $G$ acts transitively on $V_m$, for every vertex $t_1 t_2 \cdots t_m \in V_m$ there exists $g \in G$ such that $g \cdot (w_1 w_2 \cdots w_m) = t_1 t_2 \cdots t_m$. Thus the image $g \cdot w \in U_m(t_1 \cdots t_m)$, and the orbit of $w$ is dense in $\cP$. We obtain that $G$ acts minimally on $\cP$. 

Let $t = t_1 t_2 \cdots$ and $w = w_1 w_2 \cdots $ be two infinite sequences in $\cP$. We define a metric $D$ on $\cP$ by
  \begin{align}\label{eq-metric} D(t,w) = \frac{1}{d^m}, \textrm{ where }m = \max \{ n \mid t_n = w_n\}, \end{align}
that is, $D$ measures the length of the longest initial finite word contained in both $t$ and $w$. 
Since $G$ acts bijectively on $V_n$ for $n \geq 1$, $t$ and $w$ contain a common word of length $m$ if and only if the images $g \cdot t$ and $g \cdot w$ contain a common word of length $m$, so the action of $G$ on $\cP$ is equicontinuous with respect to the metric $D$,  where we can take $\delta = \epsilon$ for every $\epsilon >0$.

\subsubsection{Group chains}\label{subsec-groupchains} An important tool for studying the dynamics of equicontinuous actions on Cantor sets are group chains.

\begin{defn}\label{defn-groupchain}
Let $G$ be a countably generated discrete group. A nested descending sequence $ \{G_n\}_{n \geq 0} = G_0 \supset G_1 \supset G_2 \supset \cdots$, with $G_0 = G$, of finite index subgroups of $G$ is called a \emph{\bf group chain}.
\end{defn}

Any group chain $\{G_n\}_{n \geq 0}$ gives rise to a Cantor group action as in Example \ref{eq-groupchainaction}.

\begin{ex}\label{eq-groupchainaction}
{\rm
Let $\{G_n\}_{n \geq 0}$ be a group chain as in Definition \ref{defn-groupchain}. Then for every $n \geq 0$ the coset space $G/G_n$ is a finite set. Define  $V_n = G/G_n$ to be the set of vertices in a tree $T$. Let  $R_n = G_n/G_{n+1}$, then $V_n = G/G_n   \cong \prod_{1 \leq k \leq n}R_k$ are isomorphic as sets.

Inclusions of cosets induce the mappings
  \begin{align}\label{eq-inclus}\nu^{n+1}_n: G/G_{n+1} \to G/G_n : g G_{n+1} \to gG_n.\end{align}
Define the set of edges E in T by saying that a pair of vertices $[g_nG_n,g_{n+1}G_{n+1}]$ is an edge if and only if $g_{n+1}G_{n+1} \subset g_nG_n$. Then it is immediate that the inverse limit space 
  \begin{align}\label{eq-Ginfty}G_\infty = \lim_{\longleftarrow}\{G/G_{n+1} \to G/G_{n}\}  = \{(g_0G_0, g_1 G_1,\ldots) \mid \nu^{n+1}_{n}(hG_{n+1}) = gG_n\}\end{align}
can be identified with the space of paths $\cP \cong \prod_{n \geq 1}R_n$ of a rooted tree $T$ as in Section \ref{ex-treecylinder}. It follows that the inverse limit space $G_\infty$ is a Cantor set.

The left action of $G$ on coset spaces $G/G_n$ induces a natural left action of $G$ on $G_\infty$, given by the left multiplication
  \begin{align}\label{eq-leftaction}g \cdot (g_0 G_0, g_1 G_1,\ldots) = (gg_0G_0, gg_1G_0, \ldots). \end{align}
Denote by $(G_\infty, G)$ this action. Since $G_\infty$ is identified with $\cP$, the action \eqref{eq-leftaction} induces an action of $G$ on $\cP$, denoted by $(\cP,G)$. The group $G$ permutes the cosets in $G/G_n$, and acts transitively on each coset space $G/G_n$, so by Section \ref{ex-treecylinder} the action $(G_\infty,G)$ is minimal and equicontinuous. 
}
\end{ex}

Conversely, given a minimal equicontinuous action on the path space $\cP$ as in Section \ref{ex-treecylinder}, one can associate to it a group chain as in Example \ref{ex-clopens}.

\begin{ex}\label{ex-clopens}
{\rm
Let $T$ be a spherically homogeneous rooted tree with minimal and equicontinuous action of a discrete group $G$ as in Section \ref{ex-treecylinder}. In particular, the action of $G$ is transitive on each $V_n$, $n \geq 1$. Let $\cP$ be the space of infinite paths in $T$.

Let $x = (v_n)_{n \geq 0}$ be a path. Let $G_n = \{g \in G \mid g \cdot v_n = v_n\}$ be the subgroup of elements in $G$ which fix the vertex $v_n$, called the \emph{stabilizer} of $v_n$, or the \emph{isotropy subgroup} of the action of $G$ at $v_n$. Since $G$ acts transitively on the finite set $V_n$, then we have $|G: G_n| = |G/G_n| = |V_n|$, so $G_n$ has finite index in $G$. If $g \in G$ fixes $v_n$, then it fixes $v_i$ for $ 0 \leq i < n$, which implies that $G_n \subset G_{i}$ for $0 \leq i < n$. So the isotropy subgroups form a nested chain $\{G_n\}_{n \geq 0}$ of finite index subgroups of $G$. For each $n \geq 0$, the subgroup $G_n$ preserves the clopen set $U_n = U_n(v_n)$.
}
\end{ex}

Every minimal equicontinuous group action on a Cantor set is conjugate to a dynamical system associated to a group chain as in Definition \ref{defn-groupchain}, see \cite{ClarkHurder2013} and \cite[Appendix]{DHL2016} for details. Then Example \ref{eq-groupchainaction} shows that every such action is conjugate to an action on the path space of a spherically homogeneous tree $T$. A point $x \in X$ corresponds to a path in the path space $\cP$ of $T$, and a basic clopen set $U$ corresponds to the subset of paths in $\cP$ which contain a given vertex $v \in V$. Thus considering minimal equicontinuous actions $(X,G,\Phi)$ we may restrict to actions on path spaces of spherically homogeneous trees. 

For a given equicontinuous action $(X,G,\Phi)$, the choice of an associated chain $\{G_n\}_{n \geq 0}$ depends on a choice of a point $x \in X$, and on a choice of clopen sets $U_1 \supset U_2 \supset \cdots$. So the choice of a group chain $\{G_n\}_{n \geq 0}$ is not unique, and distinct group chains can define conjugate actions. It was shown in \cite{FO2002} that if two actions $(X,G,\Phi)$ and $(X',G',\Phi')$ are conjugate, and $\{G_n\}_{n \geq 0}$ is a group chain associated to $(X,G,\Phi)$, then there exists a sequence of indices $\{n_i\}_{i \geq 0}$ and a sequence of group elements $\{g_{n_i}\}_{i \geq 0}$ such that the group chain $\{g_{n_i} G_{n_i} g_{n_i}^{-1}\}_{i \geq 0}$ is associated to $(X',G',\Phi')$, see \cite[Section 4A]{DHL2017} for details. Thus to to study an action $(X,G,\Phi)$ in terms of group chains, it is sufficient to consider the chains of conjugate subgroups $\{g_n G_n g_n^{-1}\}_{n \geq 0}$. 

 \subsection{Ellis group for equicontinuous actions} \label{subsec-ellis}

The  \emph{Ellis (enveloping) semigroup} associated to a continuous group action $\Phi \colon G \to Homeo(X)$ on a topological space $X$   was introduced in the papers \cite{EllisGottschalk1960,Ellis1960}, and is treated in the books   \cite{Auslander1988,Ellis1969,Ellis2014}.  In this section we briefly recall some basic properties   of the Ellis group for a special case of equicontinuous   minimal group actions on Cantor sets.

 Let $X$ be a metric space, and $G$ be a countably generated group acting on $X$ via the homomorphism $\Phi: G \to {Homeo}(X)$. Suppose the action $(X, G, \Phi)$ is equicontinuous. Then the closure  $\overline{\Phi(G)}  \subset {Homeo}(X)$ in the \emph{uniform topology on maps} is identified with the Ellis group of the action. Each element of $\overline{\Phi(G)}$ is the limit of a sequence of maps in $\Phi(G)$, and we use the notation $(g_i)$ to denote a sequence $\{g_i \mid i \geq 1\} \subset G$ such that the sequence $\{\Phi(g_i) \mid i \geq 1\} \subset {Homeo}(X)$ converges in the uniform topology. 

Assume that the action of $G$ on $X$ is minimal, that is, for any $x \in X$ the orbit ${\Phi(G)}(x)$ is dense in $X$.  Then the orbit of the Ellis group $\overline{\Phi(G)}(x) = X$ for any $x \in X$. That is, the group $\overline{\Phi(G)}$ acts transitively on $X$. Denote the isotropy group of this action at $x$ by 
\begin{align}\label{iso-defn2}
\overline{\Phi(G)}_x = \{ (g_i) \in \overline{\Phi(G)} \mid (g_i) \cdot x = x\},
\end{align}
where $(g_i)\cdot x : = (g_i (x))$, for a homeomorphism $(g_i)$ in $ \overline{\Phi(G)}$. We then have the natural identification $X \cong \overline{\Phi(G)}/\overline{\Phi(G)}_x$ of left $G$-spaces.

Given an  equicontinuous   minimal Cantor system $(X,G,\Phi)$, the  Ellis group     $\overline{\Phi(G)}$  depends only on the image  $\Phi(G) \subset Homeo(X)$, while the isotropy group $\overline{\Phi(G)}_x$ of the action may depend on the point $x \in X$. 
Since the action of $\overline{\Phi(G)}$ is transitive on $X$, given any $y \in X$,   there is an element $(g_i) \in \overline{\Phi(G)}$ such that $(g_i) \cdot x  =y$. It follows that 
  \begin{align}\label{eq-profiniteconj}
  \overline{\Phi(G)}_y = (g_i) \cdot \overline{\Phi(G)}_x \cdot (g_i)^{-1}  \ .
  \end{align}
This tells us that  the \emph{cardinality} of the isotropy group $\overline{\Phi(G)}_x$ is independent of the point $x \in X$, and so  the Ellis group $\overline{\Phi(G)}$ and the cardinality of $\overline{\Phi(G)}_x$ are invariants of $(X,G,\Phi)$.

The definition of the Ellis group, given above, does not provide an easy way to compute it. In \cite{DHL2016}, we developed a method for computing $\overline{\Phi(G)}$ and $\overline{\Phi(G)}_x$ which uses group chains of Section \ref{subsec-groupchains}.

For every $G_n$ consider the \emph{core} of $G_n$, that is, the maximal normal subgroup of $G_n$ given by
\begin{equation}\label{eq-core}
C_n =  {\rm core}_{G}   \, G_n =  \bigcap_{g \in {G}} gG_ng^{-1}  ~ \subseteq ~ G_n \ .
\end{equation}
Since $C_n$ is normal in $G$, the quotient $G/C_n$ is a finite group, and  the collection  $\{C_n \}_{n \geq 0}$ forms a descending chain of normal subgroups of $G$. The inclusions $C_{n+1} \subset C_n$ induce surjective homomorphisms of finite groups, given by
  $$G/C_{n+1} \to G/C_n: gC_{n+1} \mapsto gC_n.$$
The inverse limit space
\begin{align} \label{cinfty-define}  
  C_{\infty}  & =   \lim_{\longleftarrow} \, \left\{\delta^{\ell+1}_{\ell} \colon   G/C_{\ell+1} \to G/C_{\ell}   \right\}  ~ \subset \prod_{\ell \geq 0} ~ G/C_{\ell}   
\end{align} 
is a profinite group. Also, since $G_{n+1} \subset G_n$ and $C_{n+1} \subset C_n$, there are well-defined homomorphisms of finite groups $\delta_n: G_{n+1}/C_{n+1} \to G_n/C_n$, and there is the inverse limit group
  $$\cD_x = \lim_{\longleftarrow}\{G_{n+1}/C_{n+1} \to G_{n}/C_{n}\},$$
called the \emph{discriminant group} of this action. 

\begin{thm}\cite{DHL2016}\label{Ellis-Ci}
The profinite group $C_\infty$ is isomorphic to the Ellis group $\overline{\Phi(G)}$ of the action $(X,G,\Phi)$, and the isotropy group $\overline{\Phi(G)}_x$ of the Ellis group action is isomorphic to $\cD_x$.
\end{thm}

As discussed in Section \ref{subsec-groupchains}, the group chain $\{G_n\}_{n \geq 0}$ depends on the choice of a point $x \in X$, and on the choice of a sequence of clopen sets $X=U_0 \supset U_1 \supset \cdots$ such that the set of elements in $G$ whose action preserves $U_n$ is $G_n$. Since the groups $C_n$ are normal, they do not depend on the choice of $x \in X$, but they may depend on the choice of the clopen sets $\{U_n\}_{n \geq 0}$. For any choice of $x$ and $\{U_n\}_{n \geq 0}$, the group $C_\infty$ is isomorphic to the Ellis group $\overline{\Phi(G)}$, so $C_\infty$ is independent of choices up to an isomorphism. One can think of $C_\infty$ as a choice of `coordinates' for the Ellis group $\overline{\Phi(G)}$.

Similarly, the discriminant group $\cD_x$ does not depend on choices up to an isomorphism. We note that, since $\cD_x$ is a closed subgroup of a compact group $C_\infty$, it can either be finite or an infinite profinite group which is topologically a Cantor set. 

The relationship between the cardinality of $\overline{\Phi(G)}_x$ and the properties of the action was studied in \cite{DHL2016,DHL2017}. Automorphisms of the Cantor group action $(G_\infty, G)$ (where $G$ acts on the left) are given by the right action of elements of $C_\infty$ on $G_\infty$. It is shown in \cite{DHL2016} that the automorphism group acts transitively on $G_\infty$ if and only if the isotropy group $\overline{\Phi(G)}_x \cong \cD_x$ is trivial. Thus non-triviality of $\cD_x$ is seen as an obstruction to the transitivity of the action of the automorphism group of $(G_\infty,G)$, and for this reason it was called the \emph{discriminant group} in \cite{DHL2016}. The article \cite{DHL2016} also contains examples of actions where the discriminant group is a finite non-trivial group, or a Cantor group.

\subsection{The asymptotic discriminant}\label{sec-asymptotic}

Let $(X,G,\Phi)$ be a group action on a Cantor set, let $x$ be a point and let $\{G_n\}_{n \geq 0}$ be an associated group chain, that is, the actions $(X,G,\Phi)$ and $(G_\infty,G)$ are conjugate. Recall from Section \ref{subsec-groupchains} that the groups $G_m$, $m \geq 0$, consists of elements in $G$ whose action preserves  $U_m$, and the restricted action $\Phi_m=\Phi|_{U_m}$ is the action of $G_m$. 

Set $X_m = U_m$, and consider a family of equicontinuous group actions $(X_m,G_m,\Phi_m)$. Then for each $m\geq 0$ we can compute the Ellis group of the action, and the isotropy group at $x$ as follows.
 
For each $n  \geq m \geq 0$, compute the maximal normal subgroup of $G_n$ in $G_m$ by 
 \begin{equation}\label{eq-Mn}
C^m_{n} =   {\rm core}_{G_m}   \, G_{n}  \equiv  \bigcap_{g \in {G_m}} gG_{n} g^{-1}   \subset G_m .
\end{equation}
Note that $C^m_{n}$ is the kernel of the action of $G_m$ on the quotient set $G_m/G_{n}$, and $C^0_{n}= C_{n}$.
Moreover,   for all $n > k \geq m \geq 0$, we have $\ds C^m_{n}  \subset C^k_{n} \subset G_{n} \subset G_k \subset G_m$, and $C^m_{n}$ is a normal subgroup of $G_k$. Define the profinite group 
    \begin{eqnarray}  
C^m_{k,\infty}   & \cong &  \lim_{\longleftarrow} \, \left\{G_{k}/C^m_{n} \to G_k/C^m_{n + 1}  \mid n \geq k  \right\} \\ & & = \{ (g_{n} C^m_{n})   \mid  n \geq k   \ , \  g_{k} \in G_{k} \ , \  g_{n +1} C^m_{n} = g_{n} C^m_{n}  \}     \label{Dinfty-definen}  \ .  
   \end{eqnarray}
Then   $C^m_{m,\infty}$ is the Ellis   group of the action $(X_m,G_m,\Phi_m)$, with an associated group chain $\{G_n\}_{n \geq m}$. In particular, $C^0_{0,\infty} = C_{\infty}$, defined by \eqref{cinfty-define} .

Since $G_k \subset G_m$, by definition  we have that $C^m_{k,\infty} \subset C^m_{m,\infty}$, and so  $C^m_{k,\infty}$ is a clopen neighborhood of the identity in $C^m_{m,\infty}$. 

The discriminant group  associated to the truncated group chain   $\{G_n\}_{n \geq m}$ is given by
\begin{eqnarray} 
\cD_x^m & = &    \lim_{\longleftarrow}\, \left \{  G_{n+1}/C^m_{n+1}  \to G_{n}/C^m_{n} \mid n \geq m\right\}    \subset C^m_{m,\infty} \label{eq-discquotients1} \\  & = &    \lim_{\longleftarrow}\, \left \{  G_{n+1}/C^m_{n+1}  \to G_{n}/C^m_{n} \mid n \geq k\right\}    \subset C^m_{k,\infty} , \label{eq-discquotients2}
\end{eqnarray}
where we have $\cD_x^m \subset C^m_{k,\infty}$ since $G_n \subset G_k$ for $k \geq n$. The last statement can be rephrased as saying that the discriminant group $\cD_x^m$ is contained in any clopen neighborhood of the identity in $C^m_{m,\infty}$.

To relate the discriminant groups $\cD_x^m$ and $\cD_x^k$ for $k \geq m$, we define the following maps.

For each $n \geq k \geq m \geq 0$, the inclusion $C^m_{n} \subset C^k_{n}$ induces surjective  group homomorphisms  
\begin{align}\label{eq-quotiensmn} 
\phi_{m,k}^{n} \colon   G_{n}/C^m_{n}  \longrightarrow  G_{n}/C^k_{n} \ ,
\end{align}
and the standard methods show that the maps in \eqref{eq-quotiensmn}  yield   surjective homomorphisms  of the clopen neighborhoods of the identity in $C^m_{m,\infty}$ onto the profinite groups $C^k_{k,\infty}$,
   \begin{align}\label{eq-chomeomorph} 
 \widehat{\phi}_{m,k} \colon   C^m_{k,\infty}  \to C^k_{k,\infty},
 \end{align} 
 which commute with the left action of $G$. Let $\cD_{m,k} \subset C^m_{k,\infty}$ denote the image of $\cD_x^m$ under the map \eqref{eq-chomeomorph}.
 It then follows from  \eqref{eq-quotiensmn}    that for $k > m \geq 0$,  there are    surjective homomorphisms, 
\begin{equation}\label{eq-discmapsnm2}
   \cD_x = \cD_x^0   ~ \stackrel{~  \widehat{\phi}_{0,m} ~ }{\longrightarrow} ~  \cD_{0,m} \cong \cD_x^m ~ \stackrel{~  \widehat{\phi}_{m,k} ~}{\longrightarrow} ~  \cD_x^k \ .
\end{equation}
 
 Thus, given an equicontinuous group action $(X,G,\Phi)$ on a Cantor set $G$, there is an associated sequence of surjective homomorphisms of discriminant groups \eqref{eq-discmapsnm2}, associated to the sequence of truncated group chains $\{G_n\}_{n \geq m}$, $m \geq 0$. 
 
 Since the computation of the discriminant groups above uses a group chain, associated to the action, and the group chain depends on choices, we must introduce an equivalence relation on chains of discriminant groups. Such a relation, called the \emph{tail equivalence}, was first introduced by the author in the joint work with Hurder \cite{HL2017}. Since the definition is quite technical and we do not use it directly in computations, we omit it from this paper and refer the reader to \cite{HL2017} for details. We can now introduce the notion of the asymptotic discriminant of a Cantor minimal action.

 \begin{defn}\label{def-asympdisc2}\cite{HL2017}
 Let $(X,G,\Phi)$ be an action of a countably generated group $G$ on a Cantor set $X$, and let $\{G_{n}\}_{n \geq 0}$ be an associated group chain. Then the \emph{\bf asymptotic discriminant} for $\{G_{n}\}_{n \geq 0}$ is the tail equivalence class $[\cD_x^m]_{\infty}$ of   
 the sequence of surjective group homomorphisms 
 \begin{equation}\label{eq-surjectiveDn}
[\cD_x^m]_{\infty} = \{\psi_{m,m+1} \colon \cD_x^m \to \cD_x^{m+1} \mid m \geq 1\}
\end{equation}
 defined by the discriminant groups $\cD_x^m$ for the restricted actions of $G_m$ on the clopen sets $X_m \subset X$. 
 
  \end{defn}

It was shown in \cite{HL2017} that the asymptotic discriminant is invariant under \emph{return equivalence} of group actions. Intuitively, two actions $(X_1,G,\Phi)$ and $(X_2,G,\Psi)$ are return equivalent, if there are clopen sets $U\subset X_1$ and $V \subset X_2$ such that the collections of local homeomorphisms of $U$ and $V$, induced by the actions, are compatible in a sense made precise in \cite{CHL2017}. This notion is analogous to the notion of Morita equivalence for groupoids. In the joint work with Hurder \cite{HL2017}, the author constructed an uncountable number of Cantor group actions of the same subgroup of ${\bf SL}(n,\mZ)$ with pairwise distinct asymptotic discriminants. These actions are not return equivalent.

\subsection{Locally quasi-analytic (LQA) actions} \label{subsec-LQA} Locally quasi-analytic and stable actions were defined in Definitions \ref{defn-lqa-rev} and \ref{defn-stable-wild-rev}.
 The relationship between the LQA property of group actions and their asymptotic discriminant was studied in \cite{HL2017}. In particular, Proposition 7.4 of \cite{HL2017} can be rephrased as follows.
 
 \begin{prop}\label{prop-sqa-asymptotic}\cite[Proposition 7.4]{HL2017}
 Let $(X,G,\Phi)$ be a minimal equicontinuous action of a countably generated group $G$ on a Cantor set $X$, and $x \in X$ be a point. Let $\{G_n\}_{n \geq 0}$ be an associated group chain, and let $\overline{\Phi(G)}$ be the Ellis group of the action. Then the asymptotic discriminant $[\cD_x^m]_\infty$ of the action is asymptotically constant if and only if the action of $\overline{\Phi(G)}$ on $X$ is LQA.
 \end{prop}

Thus an action $(X,G,\Phi)$ is stable if and only if its asymptotic discriminant is asymptotically constant, which means that there exists $m_0 \geq 0$ such that for all $k > m \geq m_0$, the group homomorphisms $\widehat{\phi}_{m,k}: \cD_x^m \to \cD_x^k$, defined in \eqref{eq-discmapsnm2}, are isomorphisms. 

\begin{defn}
A minimal equicontinuous group action $(X,G,\Phi)$ is \emph{\bf stable with discriminant group} $D$, if there exists $m_0 \geq 0$ such that for all $k > m \geq m_0$ there is an isomorphism $\cD_x^k \cong D$ of the discriminant groups. 
\end{defn}

\begin{ex}\label{ex-abelian-asymptotic}
 {\rm
 Suppose the group $G$ is abelian, and let $\{G_n\}_{n\geq 0}$ be a group chain in $G$. Then for $n \geq 1$ the group $G_n$ is normal in $G$, and so $C_n = G_n$ and the quotient space $G_n/C_n$ is a singleton. Then for all $m \geq 0$ the discriminant group $\cD_x^m$ is trivial and the maps $\psi_{m,m+1}$ are trivially isomorphisms. Thus $[\cD_x^m]_\infty$ is asymptotically constant. We conclude that every equicontinuous minimal action of an abelian group $G$ on a Cantor set is stable with trivial discriminant group.
 }
 \end{ex}
 
\begin{ex}
{\rm
 It was shown in \cite{DHL2017} that every finite group, and every separable profinite group can be realized as the discriminant group of a stable minimal equicontinuous group Cantor action.
}
 \end{ex} 
  
\subsection{Germinal groupoid and non-Hausdorff elements}\label{subsec-nonHaus}  
The LQA property for a group action  $(X,G,\Phi)$ can be related to the properties of the germinal groupoid $\cG(X,G,{\Phi})$ associated to the action. This groupoid  is fundamental for the study of the $C^*$-algebras these actions generate, as discussed for example by Renault in \cite{Renault1980,Renault2008}.   
 
 Recall that for $g_1, g_2 \in G$, we say that $\Phi(g_1)$ and $\Phi(g_2)$ are \emph{germinally equivalent} at $x \in X$ if $\Phi(g_1) (x) = \Phi(g_2)(x)$, and there exists an open neighborhood $x \in U \subset X$ such that the restrictions agree, $\Phi(g_1)|U = \Phi(g_2)|U$. We then write $\Phi(g_1) \sim_x \Phi(g_2)$. For $g \in G$, denote the equivalence class of  $\Phi(g)$ at $x$ by $[g]_x$.
 The collection of germs $\cG(X, G, \Phi) = \{ [g]_x \mid g \in G ~ , ~ x \in X\}$ is given the sheaf topology, and   forms an \emph{\'etale groupoid} modeled on $X$. 
 We recall the following result.
 \begin{prop}\cite[Proposition 2.1]{Winkel1983}\label{prop-hausdorff}
The germinal groupoid $\cG(X, G, \Phi)$ is Hausdorff at $[g]_x$  if and only if, for all $[g']_x \in \cG(X, G, \Phi)$ with $g \cdot x = g' \cdot x = y$, if there exists a sequence   $\{x_n\} \subset X$ which converges to $x$ such that $[g]_{x_n} = [g']_{x_n}$ for all $n$, then $[g]_{x} = [g']_{x}$. 
\end{prop}

For group Cantor actions, the following result was obtained in \cite{HL2018}.

\begin{prop}\label{prop-hausdorff2}\cite[Proposition 2.5]{HL2018}
If an action $(X,G,\Phi)$   is  LQA, then  $\cG(X, G, \Phi)$ is Hausdorff.     
 \end{prop}
 
 Thus if the groupoid $\cG(X, G, \Phi)$ is non-Hausdorff, then the action $(X,G,\Phi)$ is not LQA. 
 
In Proposition \ref{prop-hausdorff}, consider the composition of maps $h = g^{-1}\circ g'$. Since $g \cdot x = g' \cdot x$, then $h \cdot x = g^{-1} \circ g'\cdot x = x$. Denote by $[id]_x$ the germ of the identity map at $x \in X$. Then the statement of Proposition \ref{prop-hausdorff} reads as follows.

\begin{prop}\label{prop-hausdorff3}
The groupoid $\cG(X,G,\Phi)$ is Hausdorff if and only if for all $[h]_x \in \cG(X, G, \Phi)$ with $h \cdot x = x$, if there exists a sequence   $\{x_n\} \subset X$ which converges to $x$ such that $[h]_{x_n} = [id]_{x_n}$ for all $i$, then $[h]_{x} = [id]_{x}$. 
\end{prop}

Taking the contrapositive of this statement, we obtain that $\cG(X,G,\Phi)$ is a non-Hausdorff groupoid if and only if there exists a germ $[h]_x \in \cG(X, G, \Phi)$ with $h \cdot x = x$, and a sequence   $\{x_n\} \subset \fX$ which converges to $x$ such that $[h]_{x_n} = [id]_{x_n}$ for all $i$, and such that $[h]_{x} \ne [id]_{x}$. We call a representative $h$ of such a germ a \emph{non-Hausdorff element} of $\cG(X,G,\Phi)$. This is precisely Definition \ref{defn-nonH-rev} of a non-Hausdorff element in the Introduction. 

\subsection{LQA properties under inclusion and conjugacy}\label{subsec-conjugacy}
For our main theorems, we need two technical results proved in this section. Below we denote by $id$ the identity map on a Cantor set $X$.
 
 \begin{lemma}\label{lemma-inclusionwild}
 Let $(X,G,\Phi)$ be a minimal action of a (countable or profinite) group $G$ on a Cantor set $X$, and let $H$ be a subgroup of $G$.  If the action of $H$ on $X$ is not LQA then the action of $G$ on $X$ is not LQA.
\end{lemma}

\proof Since the action of $H$ is not LQA, then there exists a descending chain of clopen neighborhoods $\{U_n\}_{n \geq 0}$, such that $\cap_n U_n$ a singleton, and, for each $U_n$, an element $g_n$ such that the restriction $g_n|U_n \ne id $ and $g_n|U_{n+1} = id$. Since $H$ is a subgroup of $G$, every such $g_n$ is also in $G$, and so the action of $G$ on $X$ is not LQA.\endproof

\begin{prop}\label{prop-conjugacy}
Let $\Phi_i(G_i) \subset Homeo(X)$, $i =1,2$, be minimal equicontinuous actions of two countably generated groups on a Cantor set $X$. Suppose the closures of these actions are conjugate, that is, $\overline{\Phi_2(G_2)} = w \circ \overline{\Phi_1(G_1)} \circ  w^{-1} $, where $w \in Homeo(X)$. Then the following is true:
\begin{enumerate}
\item The action of  $\overline{\Phi_1(G_1)} $ on $X$ is LQA if and only if the action of  $\overline{\Phi_2(G_2)}$ on $X$ is LQA.
\item There is a non-Hausdorff element $g \in \overline{\Phi_1(G_1)} $ if and only if there is a non-Hausdorff element $h \in \overline{\Phi_2(G_2)} $.
\item The discriminant group $\overline{\Phi_1(G_1)}_x $ of the action of $\overline{\Phi_1(G_1)}$ on  $X$ is finite if and only if the discriminant group $\overline{\Phi_2(G_2)}_{w(x)} $ of the action of $\overline{\Phi_2(G_2)} $ on $X$ is finite. In this case both actions are stable with finite discriminant group.
\end{enumerate}
\end{prop} 
 
 \proof For $(1)$, we prove the contrapositive. Suppose $(X,G_1,\Phi_1)$ is not LQA, then there exists a collection of open sets $\{U_n\}_{n \geq 0}$ with $\bigcap U_n = \{x\}$, and a collection of elements $g_n \in \overline{\Phi_1(G_1)}$ such that $g_n|U_{n} \ne id$, while $g_n|_{U_{n+1}} = id$. Let $y = w(x)$, and $W_n = w(U_n)$, then $\{W_n\}_{n \geq 0}$ is a collection of open neighborhoods of $y$ with $\bigcap W_n = \{y\}$. Let $h_n = w g_n w^{-1}$, then $h_n|W_n \ne id$ and $h|W_{n+1} = id$. Thus $(X,G_2,\Phi_2)$ is not LQA. The converse is obtained by reversing the arrows in this argument, and the proof of $(2)$ is similar. In $(3)$, since the groups $\overline{\Phi_1(G_1)} $ and $\overline{\Phi_2(G_2)} $ of homeomorphisms are conjugate, the isotropy subgroups of their actions at $x$ and $y =w(x)$ respectively are conjugate. Thus  $\overline{\Phi_1(G_1)}_x$ is finite if and only if $\overline{\Phi_2(G_2)}_{w(x)}$ is finite. Since the discriminant groups are finite, and the homomorphisms \eqref{eq-surjectiveDn} are surjective, then the asymptotic discriminants of both actions are asymptotically constant with finite discriminant groups.
\endproof

We note that the conjugacy of the closures $\overline{\Phi_i(G_i)}$, $i =1,2$, need not imply the conjugacy of the actions of $G_i$, $i =1,2$. For instance, Nekrashevych \cite{Nekr2007} gives examples of actions of discrete non-isomorphic subgroups of $Aut(T)$ whose closures in $Aut(T)$ are conjugate. 

 \section{Wreath products and self-similar actions}\label{sec-wreathproduct}
 
In this section we recall the background on wreath products and automatic groups, which is necessary for the rest of the paper. The main reference here is Nekrashevych \cite{Nekr}. A nice concise exposition of the parts of the theory needed for work with Galois groups can be found in \cite{Jones2015}.

\subsection{Wreath products and the automorphism group of a tree}
Let $T$ be a $d$-ary rooted tree as in Section \ref{ex-treecylinder}, that is,  for $n \geq 1$, $R = R_n$ is a set with $d$ elements. Denote by $R^n = R \times  \cdots \times R$ the $n$-fold product. The set of vertices $V_n$ contains $d^n$ elements, every vertex $v \in V_n$ is connected by edges to precisely $d$ vertices in $V_{n+1}$, and every vertex in $V_{n+1}$ is connected by an edge to precisely one vertex in $V_n$. Denote by $Aut(T)$ the automorphism group of $T$. The elements of $Aut(T)$ act on the vertex sets $V_n$ by permutations in such a way that the connectedness of the tree is preserved, that is, for each $g \in Aut(T)$ and every pair $v \in V_n$ and $w \in V_{n+1}$ there is an edge $[v,w] \in E$ if and only if there is an edge $[g \cdot v, g \cdot w] \in E$. Thus $Aut(T)$ acts by homeomorphisms on the path space $\cP_d$ of $T$. For completeness, we briefly recall how to compute $Aut(T)$ from \cite{BOERT1996}. 

As in Section \ref{ex-treecylinder}, for $n \geq 1$ let ${\rm pr}_{n-1}:R^n \to R^{n-1}$ be the projection on the first $n-1$ factors in $R^n$, and let $b_{n}:  R^n \to V_n$ be a bijection such that 
  ${\rm pr}_{n-1}(b_n^{-1}(v)) = {\rm pr}_{n-1}(b_n^{-1}(w))$ if and only if $v$ and $w$ are connected by edges to the same vertex in $V_{n-1}$. Denote by $S_{d}$ the symmetric group on $d$ elements. Denote by $T_n$ the connected subtree of $T$ with the vertex set $V_0 \sqcup V_1 \sqcup \cdots \sqcup V_n$. The computation is by induction.

Note that $V_1 \cong R$, and so $Aut(T_1) = S_d$. Suppose $Aut(T_n)$ is known. Denote by $f: V_n \to S_d$ a function which assigns a permutation of $R$ to each $v \in V_{n}$, and let
  $S_{d}^{|V_{n}|} = \{f: V_{n} \to S_{d}\}$ be the set of all such functions. Then the wreath product
 \begin{eqnarray} \label{wp-formula}S_{d}^{|V_{n}|}  \rtimes Aut(T_{n})  \end{eqnarray}
  acts on $V_{n} \times R $ by
    \begin{eqnarray}\label{wp-action} (f,s)(v_{n},w) = (s(v_{n}),f(s(v_{n})) \cdot w). \end{eqnarray}
 That is, the action \eqref{wp-action} permutes the copies of $R$ in the product $V_{n} \times R$, while permuting elements within each copy of $R$ independently. Since $V_{n+1} \cong V_n \times R  \cong R^{n} \times R $, it follows that there are isomorphisms
   \begin{eqnarray}\label{eq-auttn} Aut(T_{n+1}) \cong S_{d}^{|V_{n}|}  \rtimes Aut(T_{n})  \cong   S_{d}^{|V_{n}|}  \rtimes \cdots \rtimes S_{d}^{|V_{1}|}  \rtimes S_{d}\end{eqnarray}
 of the group $Aut(T_{n+1})$ to the $(n+1)$-fold product $[S_d]^{n+1}$ of the symmetric groups $S_d$. 
 
 Next, note that there are natural epimorphisms $Aut(T_{n+1}) \to Aut(T_{n})$, induced by the projection on the second component in \eqref{wp-formula}. Thus the automorphism group of a $d$-ary rooted tree $T$ is the profinite group
 \begin{align}\label{eq-invlimitwp} Aut(T) = \lim_{\longleftarrow} \{Aut(T_{n+1}) \to Aut(T_{n}), n\geq 0\} \cong   \cdots \rtimes S_{d}^{|V_{n}|} \rtimes \cdots \rtimes S_{d}^{|V_{1}|}  \rtimes S_{d}. \end{align}
We will also denote such infinite wreath product of symmetric groups by $[S_d]^\infty$. 

\subsection{Self-similarity} The automorphism group $Aut(T)$ of a $d$-ary tree and some of its subgroups have an interesting property called self-similarity. We first describe this property for $Aut(T)$, and then we define it for the subgroups of $Aut(T)$.

Recall from Section \ref{ex-treecylinder} that the bijections $b_n: R^n \to V_n$ assign to each $v_n$ a label $b_n^{-1}(v_n) = t_1t_2 \cdots t_n$, where $t_i \in R$ for $1 \leq i \leq n$. We suppress the notation for $b_n^{-1}$ and write $v_n = t_1 t_2\cdots t_n$. Labels are assigned in such a way that $v_{n+1} = t_1 t_2 \cdots t_{n} t_{n+1}$ is joined by an edge to $v_n \in V_n$ if and only if $v_n = t_1 t_2 \cdots t_n$. Thus an infinite sequence $t = t_1 t_2 \cdots t_n \cdots$ corresponds to an infinite path in $\cP_d$. 

Let $v= v_1v_2 \cdots v_m$ be a finite word of length $m$, and denote by $vT$ a subtree of $T$ containing all paths through the vertex $v \in V_m$, that is, all paths which start with the finite subword $v$. The path space of $vT$ is a clopen subset $U_m(v)$ of $\cP_d$. Every vertex of $vT \cap V_n$ for $n \geq m$ has a label of the form $vw$, where $w$ is a word of length $n-m$ in the alphabet $R$. Every letter in $v$  or $w$ is a symbol in $R$, so there is a bijection on the sets of vertices
  \begin{align}\label{eq-treehomeo} \pi_{v}: vT\cap V  \to V: vw \mapsto w, \end{align}
which induces a homeomorphism of path spaces $\overline{\pi_v}: U_m(v) \to \cP_d$.
  
Now let $g \in Aut(T)$, and suppose $g$ maps $v \in V_m$ to a vertex $g(v) \in V_m$. The action of $g$ induces a homeomorphism $\Phi(g): \cP_d \to \cP_d$, and in particular maps the clopen set $U_m(v)$ homeomorphically onto the clopen set $U_m(g(v))$.  More precisely, for each vertex $vw \in vT$ there is a unique vertex $g(vw) \in g(v)T$, which is labelled by a word $g(v)w'$ for some finite word $w'$. Composing the bijections \eqref{eq-treehomeo} for $v$ and $g(v)$, we can define the bijection, called the \emph{section} of $g$ at $v$
  \begin{align}\label{eq-restriction-g}g|_v = \pi_{g(v)} \circ g \circ \pi_v^{-1}: V \to V: w \mapsto w', \end{align}
which defines an automorphism of the tree $T$, and so induces a homeomorphism $\Phi(g|_v):\cP_d \to \cP_d$.
  


The following definition of self-similar actions is adapted to the action of subgroups of $Aut(T)$. 

 \begin{defn}\label{defn-selfsim}\cite[Definition 1.5.3]{Nekr}
 Let $G$ be a subgroup of $Aut(T)$. Then $G$ is \emph{\bf self-similar} if for every $g \in G$ and every vertex $v$ in $T$ the map $g|_v$ defined by \eqref{eq-restriction-g} is in $G$.
 \end{defn}

Clearly $Aut(T)$ itself is self-similar. For self-similar subgroups of $Aut(T)$, we have the following representation.

Suppose $G \subset Aut(T)$ is self-similar, and let $g \in G$. Recall that $V_1$ is a set with $d$ vertices. Set $\sigma_g = g|V_1$, that is, $\sigma_g$ is a permutation of vertices in $V_1$ induced by the action of $g$. For every $v \in V_1$ we have $g|_v \in G$, so we can define a function $f_g: V_1 \to G^{|V_1|}: v \mapsto g|_{\sigma_g^{-1}(v)}$. Then $g$ acts on $T$ as an element $(f_g, \sigma_g)$ of the semi-direct product $G^{|V_1|} \rtimes S_d$, where $S_d$ denotes the symmetric group on $d$ elements. More precisely, by formula \eqref{wp-action}, if $w = (w_1 w_2 \cdots) \in \cP_d$, then $g$ acts on $w_1$ as $\sigma_g$, and on the infinite sequence $w_2 w_3 \cdots$ as $f_g(\sigma_g(w_1)) = g|_{w_1}$. Thus we can represent $g$ as a composition
   \begin{align}\label{eq-confusingnotation}g = (g|_{\sigma^{-1}_g(0)}, g|_{\sigma^{-1}_g(1)}, \ldots, g|_{\sigma^{-1}_g(d-1)}) \circ \sigma_g, \end{align}
 where $\sigma_g = (1, \sigma_g) \in G^{|V_1|} \rtimes S_d$ and $(g|_{\sigma^{-1}_g(0)}, g|_{\sigma^{-1}_g(1)}, \ldots, g|_{\sigma^{-1}_g(d-1)}) \in G^{|V_1|}$. Here $1$ denotes the trivial function in $G^{|V_1|}$ which assigns to each $v \in V_1$ the identity map of $vT$. Computing \eqref{eq-confusingnotation} we first apply the permutation $\sigma_g$ to $V_1$, and then the maps $g|_{\sigma^{-1}_g(v)}$ to the subtrees $vT$, $0 \leq v \leq d-1$. 
 
 Alternatively, we can also write $g$ as the following composition
    \begin{align}\label{eq-moreconfusingnotation}g = \sigma_g \circ  (g|_{0}, g|_{1}, \ldots, g|_{d-1}), \end{align}
 that is, when computing the action of $g$ we first apply the maps $g_v$ to the subtrees $vT$, $0 \leq v \leq d-1$, and then we apply the permutation $\sigma_g$ of $V_1$. Different sources in the literature use one or the other of these two ways to write an automorphism $g \in Aut(T)$ as a composition of two maps. In particular, \cite{Nekr,Jones2015,BFHJY} use \eqref{eq-moreconfusingnotation}, and \cite{Pink2013} uses \eqref{eq-confusingnotation}. 
 We will follow \cite{Pink2013} and mostly use \eqref{eq-confusingnotation}, as our results rely on those of \cite{Pink2013} and this notation is consistent with our definition of the wreath product in \eqref{wp-action}. Formulas \eqref{eq-confusingnotation} and \eqref{eq-moreconfusingnotation} together give the relation
  \begin{align}\label{eq-switch}  \sigma_g \circ  (g|_{0}, g|_{1}, \ldots, g|_{d-1}) = (g|_{\sigma^{-1}_g(0)}, g|_{\sigma^{-1}_g(1)}, \ldots, g|_{\sigma^{-1}_g(d-1)}) \circ \sigma_g, \end{align}
 which one can use to change from one notation to another one. 
 
Using \eqref{eq-confusingnotation}, we can write the elements of $G$ recursively, see Example \ref{ex-odometers}. 
  
\begin{ex}\label{ex-odometers}
{\rm
Let $\sigma \in S_d$ be a permutation. Then the element
  $$a = (a,1,\ldots,1)\sigma$$
first acts as $\sigma$ on the set $V_1$, and then as $a$ on the subtree $0T$, and as the identity map on every subtree $sT$, where $s \in \{1,2,\ldots, d-1\}$. This means that we have to apply $\sigma$ to the set $0T \cap V_2$, then $a$ to $00T$ and the identity map to $0sT$, $1 \leq s \leq d-1$, and then continue inductively, unfolding the action of $a$ on each next vertex level $V_n$, $ n\geq 3$. 

If $\sigma$ is a trivial permutation, then $a = (a,1,\ldots, 1)\sigma$ is the trivial map.

Next, suppose $a$ is an element in $Aut(T)$ which generates the action of an \emph{odometer}, or the \emph{adding machine}, on $\cP_d$. More precisely, start with a word $w = w_1w_2 \cdots w_n \cdots$, then
  \begin{align} \label{eq-odometer} a=\left\{ \begin{array}{ll} a(w) = (w_1+1) w_2 \cdots, & \textrm{ if }w_1 \ne d-1, \\ a(w) = 0 \cdots 0 (w_n+1) w_{n+1} \cdots, & \textrm{ if } w_i = d-1 \textrm{ for }1 \leq i \leq n-1, \textrm{ and } w_n \ne d-1, \\ a(w) = 000 \cdots, & \textrm{ if } w_n = d-1 \textrm{ for }n \geq 1.\end{array} \right. \end{align}

Clearly, $a$ acts on $V_1$ as a transitive permutation $\sigma = (0 \, 1 \cdots (d-1))$. For $v \in \{0,1,\ldots,d-2\}$ we have $a|_v = id$, and $a|_{d-1} = a$. Then, using \eqref{eq-confusingnotation},
  $$a = (a_{\sigma^{-1}(0)},a_{\sigma^{-1}(1)},\cdots, a_{\sigma^{-1}(d-1)}) \sigma = (a,1,1\cdots,1)\sigma.$$

It is convenient to represent the action of an element $g \in Aut(T)$ in a diagram, as in Figure \ref{fig-odom}. In Figure \ref{fig-odom}, $d=2$, so $T$ is a binary tree. We think of the labeling as increasing from left to right, that is, if $w_0,w_1$ are vertices in $V_{n+1}$ connected to $v \in V_n$, with $w_0$ on the left and $w_1$ on the right in the picture, then $w_1 = v0$ and $w_2 = v1$. An arc joining two edges which start at the same vertex $v$, indicates that the restriction $a|_{\sigma_g^{-1}(v)}$ is a non-trivial permutation of the set of two elements. Two shorter edges emanating from a vertex $v$ show that $a|_{\sigma_g^{-1}(v)}$ is the identity map. Even though for some vertices $a|_{\sigma_g^{-1}(v)} = 1$, this does not mean that $a$ fixes the subtree $vT$ since $a$ acts non-trivially on $V_1$ and every level $V_n$, $n \geq 1$.
Diagrams as in Figure \ref{fig-odom} are called \emph{portraits} in \cite{Nekr}.

\begin{figure}
\centering
\includegraphics[width=6cm]{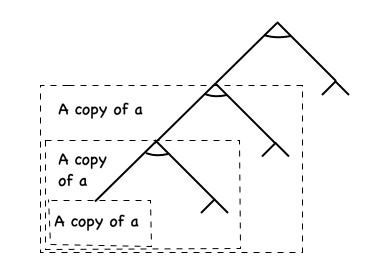}
\caption{Recursive construction of an element  $a=(a,1)\sigma$.}
\label{fig-odom}
\end{figure}

Applying \eqref{eq-switch} to $a = (a,1)\sigma$ we can write $a = \sigma(1,a)$. Then the portrait for $a$ is the mirror image of the one in Figure \ref{fig-odom}. Thus the portrait of an element $g \in Aut(T)$ depends on the representation of $g$ as a composition in \eqref{eq-confusingnotation} or \eqref{eq-moreconfusingnotation}.
}
\end{ex} 
 
Sections \eqref{eq-restriction-g} satisfy, for any finite words $v,w$, the relations \cite{Nekr}
\begin{align}\label{eq-relationsrestr} g|_{vw} = g|_v|_w, \textrm{ and } g(vw) = g(v) g|_{v}(w). \end{align}

\begin{defn}\label{defn-contracting}
Let $G \subset Aut(T)$ be a self-similar subgroup. Then $G \subset Aut(T)$ is \emph{\bf contracting}, if there is a finite set $\cN \subset G$ such that for every $g \in G$ there is $n_g \geq 0$ such that for all finite words $v$ of length at least $n_g$ we have $g|_v \in \cN$. 
\end{defn} 

The set $\cN$ is called the \emph{nucleus} of the group $G$, if $\cN$ is the smallest possible set satisfying Definition \ref{defn-contracting}. Iterated monodromy groups of post-critically finite polynomials which are the object of the study in this paper are known to be contracting \cite[Theorem 6.4.4]{Nekr}.

We also consider the following special subsets of the group $G$, introduced in \cite{Jones2015}. Let
 \begin{align}\label{eq-N0} \cN_0 = \{g \in G \mid g|_v = g \textrm{ for some non-empty }v \in V \}. \end{align}
The set $\cN_0$ is always non-empty, as it contains the identity of $G$.  It is proved in \cite[Proposition 3.5]{Jones2015} that if $G$ is contracting, then $\cN_0$ is finite and the nucleus of $G$ is given by
  $$\cN = \{h \in G \mid  h = g|_v \textrm{ for some }g \in \cN_0 \textrm{ and }v \in V \}.$$
Also define
  \begin{align}\label{eq-N1} \cN_1 = \{g \in G \mid g|_v = g \textrm{ and }g(v) = v \textrm{ for a non-empty  word } v \in V\}. \end{align} 
Then $\cN_1$ contains elements in $\cN_0$ which fix at least one path in $T$, so $\cN_1 \subset \cN_0$ and $\cN_1$ is finite for contracting actions. Also, $\cN_1$ is non-empty as it contains the identity of $G$. The following statement is proved in the last paragraph of \cite[Section 4]{Jones2015} on p. 2033.

\begin{lemma}\cite{Jones2015}\label{lemma-N1torsion}
Let $G \subset Aut(T)$ be contracting. Then every $g \in \cN_1$ is torsion.
\end{lemma}

\subsection{Non-Hausdorff elements and contracting actions}
In this section, we prove some technical results about non-Hausdorff elements in contracting subgroups of $Aut(T)$. 

Let $G \subset Aut(T)$ be a finitely generated subgroup acting on a Cantor set $\cP_d$ of infinite paths in a $d$-ary tree $T$. Here we can assume that $G$ is finitely generated, as the discrete iterated monodromy group ${\rm IMG}(f)$ associated to a post-critically finite polynomial $f(x)$ over $\mC$ is always finitely generated \cite{Nekr}.

\begin{lemma}\label{non-hausdorff-n1}
Let $G \subset Aut(T)$ be contracting, and suppose $G$ contains a non-Hausdorff element $h$. Then there is a non-Hausdorff element $g$ in $\mathcal{N}_1$.
\end{lemma}

\proof By Definition \ref{defn-nonH-rev} if $h \in G$ is non-Hausdorff, then there exists a path $\overline{x} = x_1x_2 \cdots \in \cP_d$ and a collection $\{U_n\}_{n \geq 0}$ of decreasing open neighborhoods of $\overline{x}$, with $\bigcap U_n = \{\overline{x}\}$, such that $h ( \overline{x}) = \overline{x}$, and for each $n \geq 0$ an open set $U_n$ contains an open subset $W_n$ such that $h|{W_n} = id$, but $h|{U_n} \ne id$, where $h|W_n$ and $h|U_n$ are restrictions of maps to open sets. Note that this condition implies that $w \notin W_n$, since no neighborhood of $\overline{x}$ is fixed by $h$. In particular, $W_n$ is properly contained in $U_n$.

By Section \ref{ex-treecylinder} without loss of generality we can take $U_n = U_{i_n}(\overline{x}_n)$, where $\overline{x}_n = x_1 \cdots x_{i_n}$ is a finite word labelling a vertex in the vertex set $V_{i_n}$. Then $U_n$ is the set of all paths containing the vertex $\overline{x}_n$. Also, we can take $W_n = U_{k_n}(\overline{w}_{n})$, where $\overline{w}_{n} = w_1 \cdots w_{k_n}$ is a finite word which labels a vertex in $V_{k_n}$. Since $W_n$ is properly contained in $U_n$, then we have $k_n > i_n$, and $w_j = x_j$ for $1 \leq j \leq i_n$, that is, $\overline{w}_{n} = \overline{x}_n w_{i_n+1} \cdots w_{k_n}$.

Since $h|{W_n} = id$, then $h$ fixes the word $\overline{w}_{n}$, and so it fixes the word $\overline{x}_n$. Note that the map $h|_{\overline{x}_n}$ defined by \eqref{eq-restriction-g} is a non-Hausdorff element of $G$. Indeed, set $\overline{y} = x_{i_n+1} x_{i_n+2} \cdots$, that is, $\overline{y}$ is obtained from $\overline{x}$ by discarding the first $i_n$ letters. Since $h(\overline{x}) = \overline{x}$, then $h|_{\overline{x}_n}(\overline{y}) = \overline{y}$, and so $h|_{\overline{x}_n}$ fixes every finite subword of $\overline{y}$.  The clopen sets $U_m' = U_{i_m-i_n}(x_{i_n+1} x_{i_n+2} \cdots x_{i_m})$ for $m \geq 1$ form a descending system of open neighborhoods of $\overline{y}$. Each such set contains a subset $W'_m =  U_{k_m-i_n}(x_{i_n+1} x_{i_n+2} \cdots x_{i_m}  w_{i_m+1} \cdots w_{k_m})$ fixed by $h|_{\overline{x}_n}$, while $h|_{\overline{x}_n}$ acts non-trivially on $U_m'$.

Now consider the collection of elements $\{h|_{\overline{x}_n}\}$, for $n \geq 1$. If for some $n >0$ we have $h|_{\overline{x}_n} = h$, then $h \in \cN_1$ and we are done. If not, recall that $G$ is contracting, and so by Definition \ref{defn-contracting} there is a number $\ell_h \geq 1$ such that for all finite words $v$ of length at least $\ell_h$ the restriction $h|_v \in \cN$. The collection $\{h|_{\overline{x}_n}\}$ is infinite, while the nucleus $\cN$ is finite, so there exist indices $s,t \geq 1$ such that $i_s > i_t > \ell_h$ and $h|_{\overline{x}_s} = h|_{\overline{x}_t}$. Then $g = h|_{\overline{x}_t}$ is a non-Hausdorff element in $\cN_1$.  
\endproof

Note that if $h \in \cN_1$ is non-Hausdorff, then necessarily $h \ne id$, as $h$ acts non-trivially on neighborhoods of a point $\overline{x} \in \cP_d$.

\section{The asymptotic discriminant for Chebyshev polynomials}\label{sec-cheb}

In this section, we prove Theorem \ref{thm-cheb}, which computes the asymptotic discriminant for the action of the discrete iterated monodromy group associated to a Chebyshev polynomial $T_d$ of degree $d \geq 2$, see Section \ref{sec-IMG} for a definition of a discrete iterated monodromy group.

A degree $d$ Chebyshev polynomial is defined by $T_d(x) = \cos (d \arccos x)$ for $x \in [-1,1]$, or $T_d(\theta) = \cos (d \, \theta)$ for $\theta \in [0,\pi]$. 
 Using the standard trigonometric identity for the sum of cosines, one obtains that the Chebyshev polynomials satisfy the recursive relations
   \begin{align}\label{eq-chebrecurs}T_0(x) = 1, \, T_1(x) = x, \, T_d(x) = 2xT_{d-1}(x) - T_{d-2}(x). \end{align}
  Sometimes in the literature Chebyshev polynomials are defined as
     \begin{align}\label{eq-chebrecurstilde}\widetilde{T}_0(x) = 2, \, \widetilde{T}_1(x) = x, \, \widetilde{T}_d(x) = x\widetilde{T}_{d-1}(x) - \widetilde{T}_{d-2}(x). \end{align}
The polynomials \eqref{eq-chebrecurs} and \eqref{eq-chebrecurstilde} are conjugate by a linear map, namely, $T_d(z) = \frac{1}{2}\widetilde{T}_d(2x)$.

  By \cite[Proposition 6.12.6]{Nekr}, the discrete iterated monodromy group ${\rm IMG}(T_d)$ is generated by the following maps:
  \begin{enumerate}
  \item If $d$ is even, then the generators are 
    \begin{align}\label{cheb-even-gen}a = \tau, \quad b=\sigma(b,1,\ldots, 1,a) = (b,1,\ldots,1,a)\sigma, \end{align} 
 where $\tau=(0,1)(2,3)\cdots (d-2,d-1)$ and $\sigma = (1,2)(3,4)\ldots(d-3,d-2)$. The last equality for $b$ holds since $\sigma$ fixes $0$ and $d-1$.
 \item  If $d$ is odd, then the generators are 
    \begin{align}\label{cheb-odd-gen}a = \tau(1,\ldots,1,a)= (1,\ldots,1,a)\tau, \quad b=\sigma(b,1,\ldots, 1)= (b,1,\ldots,1)\sigma, \end{align} 
 where $\tau=(0,1)(2,3)\cdots (d-3,d-2)$ and $\sigma = (1,2)(3,4)\ldots(d-2,d-1)$. The last equality for $a$ holds since $\tau$ fixes $d-1$, and the last equality for $b$ holds since $\sigma$ fixes $0$.
  \end{enumerate}
 Both $\tau$ and $\sigma$ have order $2$, so $a$ and $b$ have order $2$. The composition $\alpha = ba$ has infinite order. Indeed, if $d$ is even, then $\sigma \tau = (0,2, 4, \ldots,d-2, d-1,d-3, \ldots,3,1)$ is a $d$-cycle, and  
   $$\alpha = ba = (b,1,\ldots,1,a)\sigma\tau.$$
Denote by $e$ the identity element in ${\rm IMG}(T_d)$. Computing inductively $\alpha^k$ using formula \eqref{eq-switch}, we obtain that if $1 \leq k < d$, then $\alpha^k$ acts on the vertices in the set $V_1$ as $(\sigma \tau)^k$, and so $\alpha^k \ne e$.  
Next, by a direct computation
   $$\alpha^d  = (ba,ab,\ldots, ba, ab)(\sigma \tau)^d = (\alpha,\alpha^{-1},\ldots, \alpha, \alpha^{-1}),$$  
 so $\alpha^d$ acts as the identity permutation on the set of vertices $V_1$, and as $\alpha$ or $\alpha^{-1}$ on every subtree $vT$ of $T$, where $v \in \{0,1,\ldots,d-1\}$. So $\alpha^d \ne e$. Inductively it follows that $\alpha$ has infinite order. A similar argument shows that $\alpha = ba$ has infinite order when $d$ is odd. Summarizing, ${\rm IMG}(T_d)$ is isomorphic to the infinite dihedral group, that is, 
    \begin{align}\label{eq-dihedral}{\rm IMG}(T_d)\cong \{a,b \mid a^2 = b^2 = e\} = \{b,\alpha \mid b \alpha b^{-1} = \alpha^{-1}, b^2 = 1\}. \end{align}
    
We now compute the asymptotic discriminant for the action of ${\rm IMG}(T_d)$, $d \geq 2$, and so prove Theorem \ref{thm-cheb}. For the convenience of the reader we re-state this theorem now.

\begin{thm}\label{thm-cheb-1}
Let $T_d$ be the Chebyshev polynomial of degree $d \geq 2$ over $\mC$. Then the action of ${\rm IMG}(T_d)$ is stable with discriminant group $\mZ/2\mZ$, the finite group of order $2$.
\end{thm}

\proof From the discussion before the theorem it follows that the generator $\alpha$ in \eqref{eq-dihedral} acts transitively on $V_n$, for $n \geq 1$. 

Denote by $0^n$ the concatenation of $n$ symbols $0$, and for $n \geq 1$, consider the vertex $x_n=0^n$ in $V_n$. By \eqref{cheb-even-gen} and \eqref{cheb-odd-gen}, both when $d$ is even and when $d$ is odd, the generator $b$ fixes $x_n$. Since $\alpha$ acts transitively on $V_n$ and $|V_n| = d^n$, then the smallest power of $\alpha$ which fixes $x_n$ (and any other point in $V_n$) is $\alpha^{d^n}$. Then the isotropy group of the action of ${\rm IMG}(T_d)$ at $x_n$ is given by
  $$G_n = \{ g \in G_{\rm geom} \mid g \cdot x_n = x_n \} = \langle b, \alpha^{d^n} \rangle .$$
Since ${\rm IMG}(T_d)$ acts transitively on $V_n$, then there is a bijection ${\rm IMG}(T_d)/G_n \to V_n$ such that $eG_n \mapsto x_n$, where $e$ is the identity in ${\rm IMG}(T_d)$, with the cosets of ${\rm IMG}(T_d)/G_n$ represented by the powers $\alpha^s$, $0 \leq s \leq d^{n}-1$. In particular, $\alpha G_n \ne \alpha^{d^n-1} G_n = \alpha^{-1}G_n$ if $n \geq 2$ for any $d \geq 2$.
  
Denote by $C_n$ the set of elements which act trivially on $V_n$, then $C_n$ is the maximal normal subgroup of $G_n$ in ${\rm IMG}(T_d)$. 

\begin{lemma}\label{eq-maxcncheb}
For $n \geq 1$ the maximal normal subgroup of $G_n=\langle \alpha^{d^n},b \rangle$ in ${\rm IMG}(T_d)$ is $C_n = \langle \alpha^{d^n} \rangle$, and the discriminant group of the action of ${\rm IMG}(T_d)$ on the path space $\cP_d$ of the tree $T$ is $\cD_{{x}} \cong \mZ/2\mZ$.
\end{lemma}

\proof
Since $\alpha^{d^n}$ is the smallest power of $\alpha$ which fixes every vertex in $V_n$, then we have $\langle \alpha^{d^n} \rangle \subseteq C_n$. Consider the action of $b$ on the cosets of ${\rm IMG}(T_d)/G_n$. We have 
  $$b \alpha^s G_n  = \alpha^{-s} b^{-1} G_n = \alpha^{-s}G_n = \alpha^{d^n-s} G_n.$$
Note that $\alpha^{d^n-s} = \alpha^s$ if and only if $d^n = 2s$. So if $d$ is even,  then $b$ fixes the cosets $eG_n$ and $\alpha^{d^n/2}G_n$, and if $d$ is odd, then $b$ fixes only $eG_n$. In both cases for $n \geq 2$ the element $b$ acts non-trivially on ${\rm IMG}(T_d)/G_n$, and so $b \notin C_n$.

Since $b^2 = 1$, every non-trivial element of ${\rm IMG}(T_d)$ which is not a power of $\alpha$ is of the form $\alpha^{t_1}b\alpha^{t_2}b \cdots \alpha^{t_m}b\alpha^{t_{m+1}}$ for some $m \geq 0$, where $t_{i} \in \mZ$. Using the relations in \eqref{eq-dihedral}, this word can be reduced to the form $\alpha^t b$, for some $t \in \mZ$. For $1 \leq t \leq d^n-1$ the elements $\alpha^tb$ are not in $C_n$. Indeed, we have 
  $$\alpha^tb(e G_n) = \alpha^tG_n \ne eG_n.$$ 
It follows that the action of $a^tb$ on ${\rm IMG}(T_d)/G_n$ is trivial if and only if $t = k d^n$ for some $k \in \mZ$, and so $C_n = \langle \alpha^{d^n} \rangle$. 
  
Then $G_n/C_n = \{e C_n, bC_n\}$. The coset inclusion maps $G_{n+1}/C_{n+1} \to G_n/C_n$ are clearly bijective. It follows that 
 $$\cD_{{x}} = \lim_{\longleftarrow}\{G_{n+1}/C_{n+1} \to G_{n}/C_{n}\} \cong \mZ/2\mZ.$$
 \endproof

To compute the asymptotic discriminant we consider the restriction of the action of ${\rm IMG}(T_d)$ to clopen subsets of $\cP_d$. For $k \geq 0$, consider a truncated chain $\{G_n\}_{n \geq k}$. We have $G_k = \langle \alpha^{d^k},b \rangle$, and $G_n = \langle \alpha^{d^n},b \rangle \subset G_k$ for $n \geq k$. The group $G_k$ fixes the word $x_k = 0^k$, and so acts on the clopen subset $U_k(x_k)$ of paths through the vertex $x_k$. All such paths are in the subtree $x_kT$ of $T$. Denote by $C^k_n$ the maximal normal subgroup of $G_n$ in $G_k$, that is, $C^k_n$ contains all elements of $G_k$ which act trivially on the coset space $G_k/G_n$.

A computation similar to the one in Lemma \ref{eq-maxcncheb}, with $\alpha^{d^k}$ instead of $\alpha$, shows that $C^k_n = \langle (\alpha^{d^k})^{d^{n-k}} \rangle = \langle \alpha^{d^n} \rangle$ is the maximal normal subgroup of $G_n$ in $G_k$, the quotient $G_n/C^k_n  = \{eC^k_n,bC^k_n\} \cong \mZ/2\mZ$, and the discriminant group of the action of $G_k$ on $U_k(x_k)$ is 
  $${\ds \cD^k_{{x}} = \lim_{\longleftarrow}\{G_{n+1}/C^k_{n+1} \to G_n/C^k_n\} \cong \mZ/2\mZ}.$$
Since $C^k_n = C_n$, then the coset maps $\psi_{1,k}^{n}:G_n/C_n \to G_n/C^k_{n}$ are the identity maps, and the induced map $\psi_{1,k}:\cD_{{x}} \to \cD^k_{{x}}$ on the inverse limits is an isomorphism. Thus the asymptotic discriminant of the action of ${\rm IMG}(T_d)$ is stable with finite discriminant group, for $d \geq 2$. This finishes the proof of Theorem \ref{thm-cheb-1}.
\endproof

\section{The asymptotic discriminant of the geometric iterated monodromy group of a quadratic polynomial}\label{sec-geom}

In this section we study the geometric iterated monodromy group for post-critically finite quadratic polynomials and prove Theorem \ref{thm-1}. The proof of statement (3) of this theorem follows from Theorem \ref{thm-cheb} in Section \ref{sec-cheb}. Statement $(1)$ follows from Lemma \ref{lemma-odometer}. Before we start the proof of Theorem \ref{thm-1}, we prove a series of propositions which are used in the proof.

So let $f(x)$ be a quadratic polynomial with coefficients in the ring of integers of a number field $K$. Recall from Section \ref{sec-IMG} that $f$ induces a map $f:\mathbb{P}^1(\mC) \to \mathbb{P}^1(\mC)$ on the Riemann sphere, which has two critical points in $\mathbb{P}^1(\mC)$. One of them is the point at infinity, denoted by $\infty$, which satisfies $f(\infty) = \infty$. We denote the other point by $c$, so $C = \{c, \infty \}$. 

As in the Introduction, we denote by $P_c = \bigcup_{n \geq 1} f^n(c) $ the orbit of $c$. We assume that $P_c = \{p_1,\ldots,p_r\}$ is finite, such that $f(c) = p_1$, and $f(p_i) = p_{i+1}$ for $1 \leq i < r$. If the orbit of $c$ is strictly periodic, then $f(p_r) = p_1$. Then, since $c$ is critical, $c$ is the only preimage of $p_1$ and so $p_r = c$ and $c \in P_c$ \cite[Section 1.9]{Pink2013}. If the orbit of $c$ is strictly pre-periodic, then there exists $1 \leq s < r$ such that $f(p_r) = p_{ s+1}$. In this case the critical point $c$ is not in $P_c$. 

We consider the action of the geometric monodromy group ${\rm Gal}_{\rm geom}(f)$ on the binary tree $T$, see Section \ref{sec-IMG} for an explanation how this action arises. Recall that $\cP_2$ denotes the space of paths in the binary tree $T$, which is a Cantor set by Section \ref{ex-treecylinder}. Here the subscript in $\cP_2$ refers to the degree of the polynomial $f(x)$. 

Denote by $\sigma$ the non-trivial permutation of a set of $2$ elements. Given a set $\cA_r = \{a_1,\ldots, a_r\}$, with elements $a_i$, $1 \leq i \leq r$, to be specified later, we denote by $\widetilde{G}_r = \langle \cA_r \rangle$ a countable subgroup of $Aut(T)$ generated by $\cA_r$,  and by ${\rm CL}(\widetilde{G}_r)$ the closure of $\widetilde{G}_r$ in $Aut(T)$. Then ${\rm CL}(\widetilde{G}_r)$ is the Ellis group of the action of $\widetilde{G}_r$, see \cite{Lukina2018} for details.
 
\subsection{Strictly periodic case}
We first consider the case when the orbit of $c$ is strictly periodic.

\begin{lemma}\label{lemma-odometer}
If the critical orbit of $f(x)$ is strictly periodic with $\#P_c = 1$, then the action of ${\rm Gal}_{\rm geom}(f)$ is stable with trivial discriminant group.
\end{lemma}
\proof
By \cite[Proposition 1.9.2]{Pink2013} ${\rm Gal}_{\rm geom}(f)$ is conjugate in $Aut(T)$ to the closure of a subgroup generated by an element $a_1 = (a_1,1)\sigma$. As explained in Example \ref{ex-odometers}, the action of $a_1$ on the path space $\cP_2$ is an odometer action. Then $\widetilde{G}_r \cong \mZ$ and is abelian. Then by Example \ref{ex-abelian-asymptotic} the action of ${\rm CL}(\widetilde{G}_r)$ is stable with trivial asymptotic discriminant. By Proposition \ref{prop-conjugacy}(3) the action of ${\rm Gal}_{\rm geom}(f)$ is stable with trivial discriminant group.
\endproof

Now suppose the orbit of the critical point $c$ consists of at least two points, that is, $\#P_c \geq 2$. The class of actions with $\#P_c = 2$ includes those where ${\rm IMG}(f)$ is the well-known and well-studied Basilica group, see \cite{Nekr} or \cite{GZ2002} and references therein. The class of actions with $\#P_c = 3$ includes those associated to the polynomials whose Julia set is the `Douady rabbit' or the `airplane' \cite{Nekr}. 

\begin{prop}\label{prop-thm11}
Let $T$ be a binary tree, and let $\cA_r= \{a_1,a_2,\ldots, a_r\}$, $r \geq 2$, be the set of elements in $Aut(T)$ given by
\begin{align}\label{gen-periodic} a_1 = (a_r,1)\sigma, \, a_i = (a_{i-1},1) \textrm{ for }2 \leq i \leq r.\end{align}
Let $\widetilde{G}_r  \subset Aut(T)$ be a discrete group generated by $\cA_r$. Then the following is true:
\begin{enumerate}
\item The action of $\widetilde{G}_r $ on the space of paths $\cP_2$ is not LQA.   
\item The group $\widetilde{G}_r $ contains no non-Hausdorff elements.
\end{enumerate}
\end{prop}

\begin{figure}
\includegraphics[width=7cm]{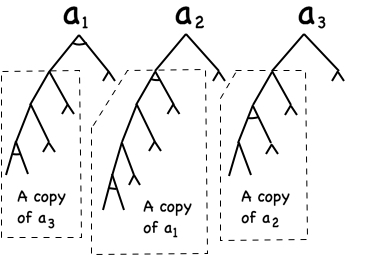}
\caption{Recursive construction of generators in the case when $r=3$ and the orbit of $c$ is strictly periodic.}
\label{fig-s0r3}
\end{figure}

\proof To show that the action of $\widetilde{G}_r $ is not LQA, we have to find a descending chain of clopen sets $\{W_n\}_{n \geq 0} \subset \cP_2$, such that ${\rm diam}(W_n) \to_n 0$, and, for each $n \geq 0$, an element $g_n \in \widetilde{G}_r $, such that the restriction $g_n|W_{n}$ is non-trivial, while the restriction $g_n|{W_{n+1}}$ is the identity map.

In the arguments below we use the labelling of vertices in $T$ by finite words of $0$'s and $1$'s as in Section \ref{ex-treecylinder}. We denote by $T_n$ the finite subtree of $T$ with vertex set $\bigsqcup_{0 \leq i \leq n} V_i$.

The generators $a_1,a_2,\ldots, a_r$ are defined recursively, so to understand how they act on the tree $T$ one has to `uncover' their action level by level. Let us start with the vertex set $V_1$ of $T$.  

By definition in formula \eqref{gen-periodic}, the generator $a_1$ acts as a non-trivial permutation (a $2$-cycle) on the vertices in $V_1$, and then it acts as $a_r$ on the infinite subtree $0T$ and as the identity map on the infinite subtree $1T$, see Figure \ref{fig-s0r3} for the portraits of the generators in the case $r=3$. So to understand how $a_1$ acts on $0T$, we have to understand how $a_r$ acts on the tree $T$. 

The generator $a_2$ acts trivially on the vertex set $V_1$, as $a_1$ on the infinite subtree $0T$, and as the identity map on the infinite subtree $1T$. Thus $a_2$ acts as a $2$-cycle on the vertices in $0T \cap V_2$, as $a_r$ on the infinite subtree $00T$ and as the identity on the infinite subtree $01T$. Denote by $0^i$ the concatenation of $i$ symbols $0$. Continuing by induction, we obtain that for $1 \leq i \leq r$, the generator $a_i$ acts as the identity on the vertex set $V_{i-1}$, as $a_1$ on the infinite subtree $0^{i-1}T$ and as the identity on every infinite subtree $wT$, where $w$ is a word of length $(i-1)$ in $0$'s and $1$'s such that at least one letter in $w$ is not $0$. 

In particular, we write 
  \begin{align}\label{eq-ar}a_r = (a_1,1,\ldots,1)1_{r-1} \end{align} 
meaning that $a_r$ acts as the identity on vertex set $V_{r-1}$, and so on the finite subtree $T_{r-1}$, as $a_1$ on the infinite subtree $0^{r-1}T$ and as the identity on every infinite subtree starting at a vertex in $V_{r-1}$ other than $0^{r-1}$. Since $a_r$ act trivially on the finite subtree $T_{r-1}$ and $a_1 = (a_r,1)\sigma$,  we have that $a_1|{T_r} = (1,1)\sigma$, that is, $a_1$ acts as $\sigma$ on $V_1$, and then as the identity on the vertices of the finite subtrees $0T \cap T_r$ and $1T\cap T_r$. It follows that for $1 \leq i \leq r$ the element $a_1$ acts on the vertex set $V_i$ as a union of $2^{i-1}$ $2$-cycles, interchanging $0$ and $1$ in the first letter of any word $w$ of length $i$, and fixing the letters from the second to the $i$-th. 

In particular, $a_1$ acts on $T_r$ as $2^{r-1}$ $2$-cycles. By \eqref{eq-switch} and \eqref{eq-ar} we obtain
\begin{align}\label{eq-a1squared} a_1^2 = (a_r,1)\sigma (a_r,1) \sigma = (a_r,1)(1,a_r) \sigma^2 = (a_r,a_r) = (a_1,1,\ldots,1,a_1, 1,\ldots,1) 1_{r}, \end{align}
so $a_1^2$ acts trivially on the first $r$ levels of the tree $T$, it acts as $a_1$ on the subtrees $0^rT$ and $10^{r-1}T$, and trivially on any other subtree $vT$, where $v$ is any word of length $r$ except $0^r$ or $10^{r-1}$.

Continuing inductively, we obtain that
  \begin{align}\label{eq-a1n}(a_1)^{2^n} = \underbrace{(a_1,1,\ldots,1, \ldots, a_1,1,\ldots, 1)}_{{\rm repeat} \, n \, {\rm times} \, (a_1,1,\ldots,1) }1_{nr}. \end{align}
Thus $(a_1)^{2^n}$ acts trivially on the first $nr$ levels of $T$. We have $|V_{nr}| = 2^{nr}$, and the pattern $(a_1,1,\ldots,1)$ of length $2^r$ is repeated $n$ times in the formula \eqref{eq-a1n}. In particular, $(a_1)^{2^n}$ acts as $a_1$ on the subtree $0^{nr}T$, and trivially on the subtree $0^{nr-1}1T$. 

So for $n \geq 1$, let $w_n = 0^{nr-1}$, and consider the set $W_n = U_{nr-1}(w_n)$ which contains all infinite sequences starting from the word $w_n$ or, alternatively,  all infinite 
paths in $\cP_2$ which pass through the vertex $w_n$. Then $W_n = U_{nr}(w_n0) \cup U_{nr}(w_n1)$, where $U_{nr}(w_n0)$ contains all paths in the subtree $0^{nr}T$, and $U_{nr}(w_n1)$ contains all paths in the subtree $0^{nr-1}1T$. By the argument above $(a_1)^{2^n}$ acts non-trivially on the clopen set $U_{nr}(w_n0)$, and trivially on the clopen set $U_{nr}(w_n1)$. 

By \cite[Proposition 2.7.1]{Pink2013} the composition $\lambda = a_1a_2 \circ \cdots \circ a_r$ generates an odometer action, and so $\lambda$ acts transitively on every level of the tree $T$. Since $|V_{nr-1}| = 2^{nr-1}$, then the power $\lambda^{2^{nr-1}}$ fixes every vertex in $V_{nr-1}$, and acts as $2^{nr-1}$ $2$-cycles on the vertices in $V_{nr}$. In particular, $\lambda^{2^{nr-1}}$ fixes $w_n$, and permutes $w_n0$ and $w_n1$. This means that $\lambda^{2^{nr-1}}$ maps $U_{nr}(w_n0)$ onto $U_{nr}(w_n1)$, and $U_{nr}(w_n1)$ onto $U_{nr}(w_n0)$. Define
  \begin{align}\label{eq-gnelement} g_n =  \lambda^{-2^{nr-1}} \circ (a_1)^{2^n} \circ \lambda^{2^{nr-1}},\end{align}
then $g_n$ acts trivially on the clopen set $U_{nr}(w_n0)$, and non-trivially on the clopen set $U_{nr}(w_n1)$, and so on $W_n$.  

Note that since $r \geq 2$, then $W_{n+1} = U_{{(n+1)}r-1}(w_{n+1}) \subset U_{nr}(w_n0) \subset W_n$, so $\{W_n\}_{n \geq 0}$ is a decreasing sequence of clopen sets, such that $g_n|W_n$ is non-trivial, while $g_n|W_{n+1}$ is trivial. We conclude that the action of $\widetilde{G}_r$ is not LQA, which proves $(1)$.

To show $(2)$ note that the group $\widetilde{G}_r$ is contracting. Then by Lemma \ref{non-hausdorff-n1}  if $g \in \widetilde{G}_r$ is non-Hausdorff then the finite set $\cN_1$ defined by \eqref{eq-N1} contains a non-Hausdorff element. By Lemma \ref{lemma-N1torsion} every element in $\cN_1$ is torsion. But the group $\widetilde{G}_r$ is a group of type $\mathfrak{R}(v)$ in \cite{BN2008}, where $v = 0^{n-1}$, and so by \cite[Proposition 3.11]{BN2008} it is torsion free. Therefore, the set $\cN_1$ for the action of $\widetilde{G}_r$ contains only the identity element, and so $\widetilde{G}_r$ does not contain any non-Hausdorff elements.
\endproof

\begin{remark}\label{remark-torsionclosure}
{\rm
Although the group $\widetilde{G}_r$ in Proposition \ref{prop-thm11} does not contain any non-Hausdorff elements, we cannot rule out that the closure ${\rm CL}(\widetilde{G}_r)$ of the action does contain them. Indeed, the absence of non-Hausdorff elements in the group $\widetilde{G}_r$ in Proposition \ref{prop-thm11} is a consequence of the fact that $\widetilde{G}_r$ is torsion free. By a celebrated result of Lubotzky \cite{Lubotzky1993} profinite completions of torsion free groups may have non-trivial torsion elements. The construction of Lubotzky was used in \cite{DHL2017} to construct examples where the action $\Phi:G \to Homeo(X)$ is that of a torsion free group, while the closure of the action $\overline{\Phi(G)}$ contains torsion elements. This motivates Problem \ref{prob-closurenonHausdorff} in the Introduction.
}
\end{remark}

\subsection{Strictly pre-periodic case.} Now suppose that the post-critical orbit $P_c = \{p_1,\ldots,p_r\}$ of the critical point $c$ of $f(x)$ is strictly pre-periodic, that is, there exists $s \geq 1$ such that $f(p_r) = p_{s+1}$.  In Proposition \ref{prop-thm12} below we consider the case when $r \geq 3$. 

Recall that we denoted by $\sigma$ the non-trivial permutation of a set of $2$ elements. Given a set $\cB_r = \{b_1,\ldots, b_r\}$, with elements $b_i$, $1 \leq i \leq r$, to be specified in Proposition \ref{prop-thm12}, we denote by $\widetilde{H}_r = \langle \cB_r \rangle$ a countable subgroup of $Aut(T)$ generated by $\cB_r$,  and by ${\rm CL}(\widetilde{H}_r)$ the closure of $\widetilde{H}_r$ in $Aut(T)$. Then ${\rm CL}(\widetilde{H}_r)$ is the Ellis group of the action of $\widetilde{H}_r$.

\begin{prop}\label{prop-thm12}
Let $T$ be a binary tree.
Let $r \geq 3$, let $1 \leq s <r$, and let $\cB_r= \{b_1,b_2,\ldots, b_r\}$ be the set of elements in $Aut(T)$ given by
\begin{align}\label{gen-preperiodic} b_1 = \sigma, \, b_{s+1} = (b_s,b_r), \, b_i = (b_{i-1},1) \textrm{ for }i \ne 1, s+1.\end{align}
Let $\widetilde{H}_r \subset Aut(T)$ be a group generated by $\cB_r$. Then $\widetilde{H}_r$ contains a non-Hausdorff element, and so the action of $\widetilde{H}_r$ on the space of paths $\cP_2$ is not LQA.
\end{prop}

\proof

We will consider two cases, first when $s+1 = r$ and so the periodic part of the orbit of the critical point $c$ is just a fixed point, and second when $r > s+1$, so that the periodic part of the orbit of $c$ has length at least $2$. 

\begin{lemma}\label{pre-periodic-fixedpoint}
In Proposition \ref{prop-thm12}, suppose that $s+1 = r$. Then $b_r$ is non-Hausdorff.
\end{lemma}

\proof If $b_r$ is non-Hausdorff, then there exists an infinite path $\overline{x} \in \cP_2$, a descending collection of clopen neighborhoods $\{W_n\}_{n \geq 1}$ with $\bigcap_{n \geq 1} W_n = \{\overline{x}\}$ and, for each $n \geq 1$, a clopen subset $O_n \subset W_n$, such that $b_r(\overline{x}) = \overline{x}$, $b_r|O_n$ is the identity, while $b_r|W_n$ is non-trivial. We will find such $\overline{x}$, $\{W_n\}_{n \geq 1}$ and $\{O_n\}_{n \geq 1}$. 
 
\begin{figure}
\includegraphics[width=12cm]{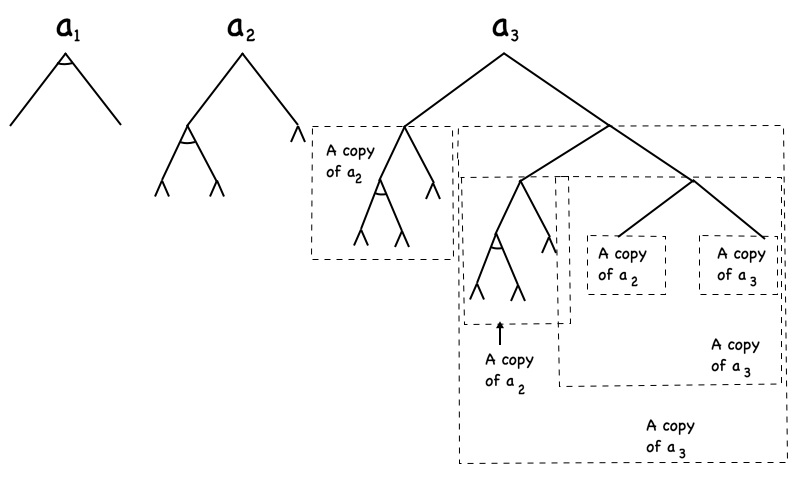}
\caption{Recursive construction of generators in pre-periodic case when $r=3$ and $s=2$.}
\label{fig-s2r3}
\end{figure}

Let us first understand how the generators $b_i$, $1 \leq i \leq s$, act on $\cP_2$. The portraits of the generators in $\cB$ for the case $s+1=r=3$ are shown in Figure \ref{fig-s2r3}. The generator $b_1 = \sigma$, so $b_1$ acts on $V_1$ as a $2$-cycle. For $n \geq 1$, $V_n$ contains $2^n$ vertices, so $b_1$ acts on $V_n$ as $2^{n-1}$ $2$-cycles. Thus $b_1$ has order $2$ and no fixed points. Note that since $r \geq 3$ and $s+1 = r$, then $s \geq 2$.

For $1 < i \leq s$ we have $b_i = (b_{i-1},1)$. That is, $b_2$ acts trivially on the vertex set $V_1$, as $b_1$ on the subtree $0T$ of $T$, and trivially on the subtree $1T$. So $b_2$ fixes a clopen set $U_1(1)$, and acts as $2^{n-2}$ $2$-cycles on the intersection $0T \cap V_n$, where $V_n$ is the vertex set at level $n \geq 2$. Inductively, one obtains that $b_i$ acts trivially on all vertices in the subtree $T_{i-1}$, as $2^{n-i}$ $2$-cycles on the vertices in the intersection $0^{i-1}T \cap V_n$, for $n \geq i$, and trivially on the rest of the tree. So $b_i$ has order $2$ for $1 < i \leq s $. In particular, $b_s$ acts non-trivially on a clopen subset of $U_1(0)$ and trivially on $U_1(1)$.

Now consider $b_r = (b_s, b_r)$, where $s +1 = r$. We will unravel how $b_r$ acts on the path space $\cP_2$ by induction on the level $n \geq 1$ in the tree $T$. From the definition, $b_r$ acts trivially on the vertex set $V_1$, as $b_s$ on $0T$ and as $b_r$ on $1T$. Then $b_r$ acts non-trivially on a clopen subset of $U_2(00)$, and trivially on $U_2(01)$. Since $b_r$ acts as $b_r$ on $1T$, then it acts as $b_s$ on $10T$ and as $b_r$ on $11T$. This means that $b_r$ acts non-trivially on a clopen subset of $U_3(100)$ and trivially on $U_3(101)$. 

Inductively, we obtain that $b_r$ acts as $b_r$ on the subtree $1^n T$, for $n \geq 1$. All infinite paths contained in this subtree are in the clopen set $W_n = U_n(1^n)$, where $1^n$ denotes a word obtained by a concatenation of $n$ copies of $1$. Then $b_r$ acts as $b_s$ on the subtree $1^n0T$, and the clopen set $U_{n+1}(1^n0)$, containing all paths of $1^n0T$. More precisely, the action of $b_r$ is non-trivial on a clopen subset of $U_{n+2}(1^n00)$, and it is trivial on the clopen subset $O_n = U_{n+2}(1^n01)$. We constructed the collections $\{W_n\}_{n \geq 1}$ and $\{O_n\}_{n \geq 1}$. 

Note that the for $n \geq 1$, we have $\bigcap_{1 \leq  k \leq n} W_k = W_n \ne \emptyset$, so $\{W_n\}_{n \geq 1}$ is a family of closed sets in $\cP_2$ with finite intersection property. Since $\cP_2$ is compact, then $\bigcap_{n \geq 1} W_n$ is non-empty \cite[Section 17]{Willard}. Any sequence $\overline{y}$ which contains at least one letter $0$ is not in $W_n$ for $n$ large enough, so it follows that $\bigcap_{n \geq 1} W_n = \overline{x} = 1^\infty$, where $1^\infty$ denotes an infinite sequence of $1$'s. By the discussion above, the action of $b_r$ is non-trivial only on subsets contained in the sets of the form $U_n(1^{n-2}00)$, for $n \geq 1$. Since $\overline{x}$ does not contain any $0$'s, it must be a fixed point of $b_r$. We have shown that $b_r$ is non-Hausdorff. Note that $b_r$ is torsion of order $2$.
\endproof

We now consider the second case, when $r > s+1$ and $s \geq 1$. 

\begin{lemma}\label{pre-periodic-orbit}
In Proposition \ref{prop-thm12}, suppose that $r> s+1$. Then for $s+1 \leq i \leq r$, the element $b_i$ is non-Hausdorff.
\end{lemma}

\proof We have that $b_1 = \sigma$, and for $1 \leq i \leq s$, the generator $b_i = (b_{i-1},1)$ acts on the tree $T$ in the same way as in Lemma \ref{pre-periodic-fixedpoint}. In particular, $b_s$ acts non-trivially on the set $U_1(0)$ if $s \geq 2$, and on the whole space $\cP_2$ if $s=1$. The generator $b_s$ has order $2$.
  
Next, we have $b_{s+1} = (b_s,b_r)$, and $b_i = (b_{i-1},1)$ for $s+2 \leq i \leq r$. Inductively, for $s+1< i  \leq r$ the generator $b_i=(b_{i-1},1)$ acts as $b_{s+1}$ on the subtree $0^{i-s-1}T$, and so on the clopen set $U_{i-s-1}(0^{i-s-1})$. Therefore, $b_r$ acts as $b_{s+1}$ on the clopen set $U_{r-s-1}(0^{r-s-1})$. Since $b_{s+1} = (b_s,b_r)$, the generator $b_r$ acts non-trivially (as $b_s$) on the clopen set $U_{r-s}(0^{r-s})$. It acts as $b_r$ on the clopen set $U_{r-s}(0^{r-s-1}1)$, and so, since $b_r = (b_{r-1},1)$, $b_r$ acts trivially on the clopen set $U_{r-s+1}(0^{r-s-1}11)$.

\begin{figure}
\includegraphics[width=12cm]{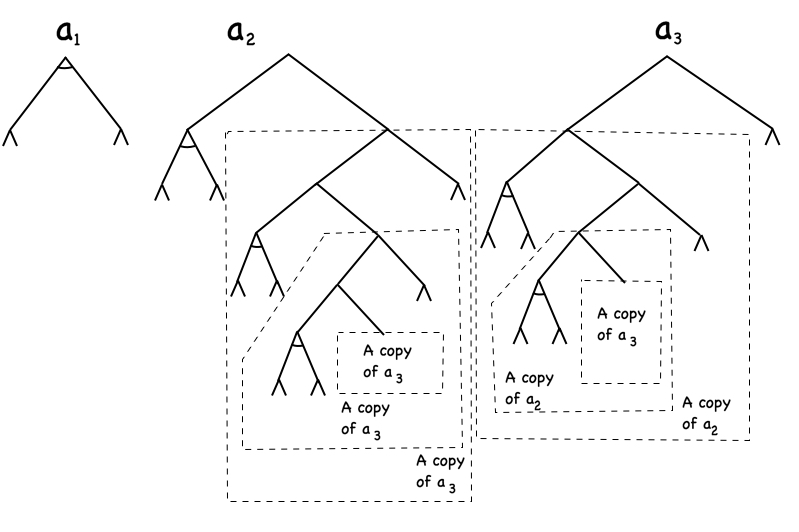}
\caption{Recursive construction of the generators in the pre-periodic case when $s=1$ and $r=3$.}
\label{fig-s1r3}
\end{figure}

Portraits of generators in the case $s=1$ and $r=3$ are presented in Figure \ref{fig-s1r3}.

For $n \geq 1$ denote by $(0^{r-s-1}1)^n$ the concatenation of $n$ copies of the word $0^{r-s-1}1$. By induction, we obtain that $b_r$ acts as $b_{s+1}$ on the clopen set $U_{n(r-s)-1}((0^{r-s-1}1)^{n-1}0^{r-s-1})$. Then it acts non-trivially (as $b_s$) on the clopen set 
  \begin{align}\label{eq-bsactionbr}Z_n = U_{n(r-s)}((0^{r-s-1}1)^{n-1}0^{r-s} \end{align}
and as $b_r$ on the clopen set 
   $$W_n = U_{n(r-s)}((0^{r-s-1}1)^n).$$
Then $b_r$ acts trivially on the clopen subset      
     $$O_n = U_{n(r-s)+1}((0^{r-s-1}1)^{n}1) \subset W_n.$$
Note that $Z_n \subset W_{n-1}$, so we obtained a nested family of clopen sets $\{W_n\}_{n \geq 1}$ such that for all $n \geq 1$, $b_r|W_n$ is non-trivial, while $b_r|O_n$, for $O_n \subset W_n$, is trivial.    
By an argument similar to the one at the end of Lemma \ref{pre-periodic-fixedpoint} we obtain that the intersection $\bigcap_{n \geq 1} W_n$ is a point $\overline{x} = (0^{r-s-1}1)^\infty$, that is, $\overline{x}$ is a concatenation of an infinite number of repetitions of the word $0^{r-s-1}1$. By the definition of $b_i$ for $s+1 < i \leq r$ one can see that $b_r$ acts non-trivially only on the subsets of the form \eqref{eq-bsactionbr}, in particular, every sequence in $Z_n$ contains a word $0^{r-s}$. Since $\overline{x}$ does not contain such a word, $\overline{x}$ is fixed by the action of $b_r$. We conclude that $b_r$ is a non-Hausdorff element.

For $s+1 \leq i <r$, by construction the path space $\cP_2$ contains a clopen neighborhood $W$ such that $b_i$ acts on $W$ as $b_{s+1} = (b_s,b_r)$. Then $W$ contains a subset $W'$ such that $b_i$ acts on $W'$ as $b_r$.  Since $b_r$ is non-Hausdorff,  then $b_i$ is non-Hausdorff.
\endproof

This finishes the proof of Proposition \ref{prop-thm12}.
\endproof

\subsection{Proof of Theorem \ref{thm-1}}

\proof Statement $(1)$ of the theorem follows immediately from Lemma \ref{lemma-odometer}, that is, if the orbit of the critical point $c$ of the polynomial $f(x)$ is a single point, then the action of ${\rm Gal}_{\rm geom}(f)$ is conjugate to the action of the enveloping group of an odometer action. By Proposition \ref{prop-conjugacy} this implies that the action of ${\rm Gal}_{\rm geom}(f)$ is stable with trivial discriminant group.

Let us prove statement $(2)$. That is, suppose the critical point $c$ of the polynomial $f(x)$ has a strictly periodic orbit of length $r = \#P_c \geq 2$. Then by \cite[Theorem 2.4.1]{Pink2013} ${\rm Gal}_{\rm geom}(f)$ is conjugate in $Aut(T)$ to the profinite group ${\rm CL}(\widetilde{G}_r)$, where $\widetilde{G}_r$ is as in Proposition \ref{prop-thm11}. Since the action of $\widetilde{G}_r $ is not LQA by Proposition \ref{prop-thm11}, then the action of the closure ${\rm CL}(\widetilde{G}_r)$ is not LQA. Then by Proposition \ref{prop-conjugacy} the action of ${\rm Gal}_{\rm geom}(f)$ on the space of paths $\cP_2$ in a tree $T$ is not LQA, and so has wild asymptotic discriminant.

For statement $(3)$, suppose the orbit of $c$ is strictly pre-periodic of cardinality $\#P_c =2$. Then by \cite[Proposition 3.4.2]{Pink2013}  ${\rm Gal}_{\rm geom}(f)$ is conjugate to the closure of the action of the group $\widetilde{H}_2$, which is generated by the set $\cB_2 =\{ b_1= \sigma, \, b_2 = (b_1,b_2)\}$. The elements in $\cB_2$ are conjugate to the generating set in \eqref{cheb-even-gen} for $d=2$. Indeed, note that $\sigma b_2 \sigma = (b_2,b_1)$ which gives the second generator in \eqref{cheb-even-gen}. Then by Proposition \ref{prop-conjugacy}(3) it follows from Theorem \ref{thm-cheb-1} that the action of ${\rm Gal}_{\rm geom}(f)$ on the space of paths $\cP_2$ in a tree $T$ is stable with finite discriminant group.

For statement $(4)$, suppose the polynomial $f(x)$ has a strictly pre-periodic orbit of cardinality $\#P_c \geq 3$. Then by \cite[Theorem 3.4.1]{Pink2013} ${\rm Gal}_{\rm geom}(f)$ is conjugate in $Aut(T)$ to the profinite group ${\rm CL}(\widetilde{H}_r)$, where $\widetilde{H}_r$ is as in Proposition \ref{prop-thm12}. Since the action of $\widetilde{H}_r$ is not LQA by Proposition \ref{prop-thm12}, then the action of the closure ${\rm CL}(\widetilde{H}_r) \subset Aut(T)$ is not LQA. Then by Proposition \ref{prop-conjugacy} the action of ${\rm Gal}_{\rm geom}(f)$ on the space of paths $\cP_2$ in a tree $T$ is not LQA, and so has wild asymptotic discriminant.
\endproof

\section{Asymptotic discriminant for the arithmetic iterated monodromy group of post-critically finite polynomials}\label{sec-arith}

Let $f(x)$ be a post-critically finite quadratic polynomial, and denote by ${\rm Gal}_{\rm geom}(f)$ and ${\rm Gal}_{\rm arith}(f)$ the geometric and the arithmetic iterated monodromy groups respectively, see Section \ref{sec-IMG} for details about these groups. In this section, we compute the asymptotic discriminant for the arithmetic iterated monodromy group, thus proving Theorem \ref{thm-3}.

Recall from Section \ref{sec-IMG} that $f(x)$ induces the map of the Riemann sphere $\mathbb{P}^1(\mC)$, which has two critical points, the point at infinity $\infty$ and the point $c \in \mC$. Recall that we denote by $P_c$ the orbit of the critical point $c$ of $f(x)$. 

\begin{prop}\label{prop-semidirect-0}
Let $K$ be a finite extension of $\mQ$, and let $f(x)$ be a quadratic polynomial with coefficients in the ring of integers of $K$, such that the orbit $P_c$ of the critical point $c$ of $f(x)$ is finite of length $r = \#P_c$. Consider the action of the profinite iterated arithmetic monodromy group ${\rm Gal}_{\rm arith}(f)$ on the space of paths $\cP_2$ of the tree $T$ of the solutions to $f^n(x) = t$, $n \geq 1$. Suppose $f(x)$ falls within one of the following categories:
\begin{enumerate}
\item either the orbit of the critical point $c$ is strictly periodic with $r \geq 2$,
\item or the orbit of the critical point $c$ is strictly pre-periodic with $r \geq 3$.
\end{enumerate} 
Then the action of ${\rm Gal}_{\rm arith}(f)$ is not LQA, and so has wild asymptotic discriminant.
\end{prop}

\proof In both cases the action of ${\rm Gal}_{\rm geom}(f)$ is not LQA by Theorem \ref{thm-1}. By \cite{Pink2013} ${\rm Gal}_{\rm geom}(f)$ is a normal subgroup of ${\rm Gal}_{\rm arith}(f)$, and both act minimally on the path space $\cP_2$ of the binary tree $T$. Then it follows by Lemma \ref{lemma-inclusionwild}  that the action of ${\rm Gal}_{\rm arith}(f)$ is not LQA.
\endproof

Suppose $f(x)$ does not satisfy Proposition \ref{prop-semidirect-0}, and the critical point $c$ of $f(x)$ has a strictly periodic orbit which is a fixed point, or when it has a strictly pre-periodic orbit of length $2$. In both cases the asymptotic discriminant of the action of ${\rm Gal}_{\rm geom}(f)$ is stable by Theorem \ref{thm-1}, so we cannot use Lemma \ref{lemma-inclusionwild} and have to compute the asymptotic discriminant of the action of ${\rm Gal}_{\rm arith}(f)$ directly. We do that below in a series of propositions.

\begin{prop}\label{prop-semidirect}
Let $K$ be a finite extension of $\mQ$, and let $f(x)$ be a quadratic polynomial with coefficients in the ring of integers of $K$. Suppose the orbit of the critical point $c$ of $f(x)$ is strictly periodic with $\#P_c = 1$. Then the action of the profinite iterated arithmetic monodromy group ${\rm Gal}_{\rm arith}(f)$ on the space of paths $\cP_2$ of the tree $T$ of the solutions to $f^n(x) = t$, $n \geq 1$ is stable with infinite discriminant group.
\end{prop}

\proof Under the hypothesis of the proposition, ${\rm Gal}_{\rm geom}(f)$ is conjugate to the closure of an odometer action, as it is explained in Lemma \ref{lemma-odometer}.  More precisely, let $a_1 = (a_1,1)\sigma$, let $\widetilde{G}_1 = \langle a_1 \rangle$ be a  countable subgroup of $Aut(T)$ generated by $a_1$, and let ${\rm CL}(\widetilde{G}_1)$ be the closure of the action of $\widetilde{G}_1$ in $Aut(T)$.  Then by  \cite[Theorem 2.8.2]{Pink2013} there is an element $w \in Aut(T)$, such that 
   \begin{align}\label{eq-wconj}{\rm CL}(\widetilde{G}_1) = w\, {\rm Gal}_{\rm geom}(f)\, w^{-1}, \end{align}
and ${\rm Gal}_{\rm geom}$ is topologically generated by an element $z$ conjugate to $w^{-1}a_1w$ in ${\rm Gal}_{\rm geom}$. That is, there is $s \in {\rm Gal}_{\rm geom}(f)$ such that $z = s^{-1}w^{-1} a_1 ws$. Then $z$ is a generator of an odometer action on the tree $T$, that is, $G_{\rm geom}= \langle z \rangle$ is a dense subgroup of ${\rm Gal}_{\rm geom}(f)$ which acts transitively on every vertex level $V_n$, $n \geq 0$, of $T$, see \cite[Proposition 1.6.2]{Pink2013}. Relabelling the vertices of the tree $T$ via the conjugating automorphism $ws$, we may assume that $z = a_1$. Then ${\rm Gal}_{\rm geom}(f) = {\rm CL}(\widetilde{G}_1) $ and $G_{\rm geom} = \widetilde{G}_1$. 

The group ${\rm Gal}_{\rm geom}(f)$ is isomorphic to the profinite group $\mbZ_2$ of the dyadic integers. 
 The normalizer $N$ of  ${\rm Gal}_{\rm geom}(f) \cong \mbZ_2$ in $Aut(T)$ is isomorphic to the semi-direct product $\mbZ_2 \rtimes \mbZ_2^\times$ \cite[Proposition 1.6.3]{Pink2013}, where $\mbZ_2^\times \cong Aut(\mbZ_2)$ denotes the multiplicative group of the dyadic integers. By \cite[Proposition 1.6.4]{Pink2013} this isomorphism can be given explicitly by
   $$\mbZ_2 \rtimes \mbZ_2^\times \to N: (m, \ell) \mapsto a_1^m z_\ell,$$
 where $m = (m_i) \in \mbZ_2$, for $\ell = (\ell_i) \in \mbZ_2^\times$ we have 
   $$z_\ell = (z_\ell, a_1^{\frac{\ell-1}{2}}z_\ell),$$ 
and the element $z_\ell$ acts on $\mbZ_2$ by raising $a_1$ to the $\ell$-th power,
      \begin{align}\label{eq-normaction1}z_\ell a_1 z_\ell^{-1} = a_1^\ell. \end{align}

Recall that  ${\rm Gal}_{\rm geom}(f)$ is a normal subgroup of ${\rm Gal}_{\rm arith}(f)$, so ${\rm Gal}_{\rm arith}(f) \subseteq N$. Recall also that we have an exact sequence \eqref{exact-1}
$$ \xymatrix{ 0 \ar[r] & {\rm Gal}_{\rm geom}(f) \ar[r] &  {\rm Gal}_{\rm arith}(f) \ar[r] & {\rm Gal}(L/K) \ar[r] & 0},$$
so ${\rm Gal}_{\rm arith}(f)/{\rm Gal}_{\rm geom}(f) \cong {\rm Gal}(L/K) $. Here the field $K$ is the base field of the polynomial $f(x)$, and $L$ is defined in Section \ref{sec-IMG}. So there is a homomorphism (see \cite[Theorem 2.6.8]{Pink2013} and \cite[Theorem 2.8.4]{Pink2013} for $r=1$)
  $$\bar{\rho}: {\rm Gal}(L/K) \to {\rm Gal}_{\rm arith}(f)/{\rm Gal}_{\rm geom}(f) \subseteq N/{\rm Gal}_{\rm geom}(f) \cong \mbZ^\times_2,$$
and elements of ${\rm Gal}_{\rm arith}(f)/{\rm Gal}_{\rm geom}(f)$ act on ${\rm Gal}_{\rm geom}(f)$ via \eqref{eq-normaction1}.

We now use the fact that $K$ is a finitely generated extension of $\mQ$. By \cite[Theorem 2.8.4]{Pink2013} in this case the homomorphism $\bar{\rho}$ is surjective onto $N/{\rm Gal}_{\rm geom}(f)$, that is, ${\rm Gal}_{\rm arith}(f) \cong \mbZ_2 \rtimes \mbZ_2^\times$ . We compute the asymptotic discriminant of the action of ${\rm Gal}_{\rm arith}(f)$ in Lemma \ref{lemma-discr0}.

In the proposition below, we take $a = a_1$, to simplify the notation. That is, the element $a$ in \eqref{eq-garith} generates the odometer action of $G_{\rm geom}$.

\begin{lemma}\label{lemma-discr0}
Let $K$ be a finite extension of $\mQ$, and let $f(x)$ be a quadratic polynomial with coefficients in the ring of integers of $K$. Suppose the orbit of the critical point $c$ of $f(x)$ is strictly periodic with $\#P_c = 1$. Then the profinite iterated arithmetic monodromy group ${\rm Gal}_{\rm arith}(f)$ is the closure of the subgroup
  \begin{align}\label{eq-garith}G_{\rm arith} \cong \langle a,b,c \mid b^2 = 1, \, bab^{-1} = a^{-1}, \, cac^{-1} = a^5, \, bcb^{-1}c^{-1} = 1\rangle \subset Aut(T), \end{align}
and the action of ${\rm Gal}_{\rm arith}(f)$ on the path space $\cP_2$ of the binary tree $T$ has an associated group chain $\{G_n\}_{n \geq 0}$, where for $n \geq 0$
  \begin{align}\label{eq-gngenerators-Q}G_n = \langle a^{2^n},b,c \rangle \subset G_{\rm arith}. \end{align}
\end{lemma}
   
 \proof As discussed just before the lemma, ${\rm Gal}_{\rm arith}(f) \cong \mbZ_2 \rtimes \mbZ_2^\times$, where $\mbZ_2$ is the dyadic integers, and $\mbZ_2^\times$ is the multiplicative group of $\mbZ_2$. We are going to determine the generators and relations for $G_{\rm arith}$. We start by building a bijection between the path space $\cP_2$ and the dyadic integers ${\ds \mbZ_2 = \lim_{\longleftarrow}\{\mZ/2^{n+1} \mZ \to \mZ/2^n \mZ\}}$.
 
 To this end, recall that vertices in $V_n$ are labelled by words of length $n$ in $0$'s and $1$'s, and let $x_n = 0^n$ and $\overline{x} = 0^\infty$. The group $G_{\rm geom} = \langle a \rangle \cong \mZ$ acts transitively on each vertex set $V_n$, $n \geq 1$, so there is a bijection $\kappa_n:V_n \to \mZ/2^n\mZ$, such that, given $y_n = t_1 t_2 \ldots t_n \in V_n$, we have  
  \begin{align}\label{eq-kappa}\kappa_n(y_n) = k +2^n \mZ \textrm{ if and only if }y_n = x_n \cdot a^k, \end{align}
where $\cdot$ denotes the action of $G_{\rm geom}$ on $V_n$, and $0 \leq k \leq 2^n-1$. In particular, $x_n$ is mapped onto the coset of $0$ in $\mZ/2^n\mZ$. It is straightforward to check that the maps $\kappa_n$, $n \geq 1$, are compatible with the bonding maps 
  $$V_{n+1} \to V_n: t_1t_2\ldots t_n t_{n+1} \mapsto t_1 t_2 \ldots t_n$$
and with the coset inclusions  $\mZ/2^{n+1}\mZ \to \mZ/2^n \mZ$. That is, for every $y_{n+1} = t_1 t_2 \ldots t_n t_{n+1} \in V_{n+1}$ we have
  $$\kappa_{n+1}(t_1 t_2 \ldots t_n t_{n+1}) \mod 2^n  = \kappa_n(t_1 t_2 \ldots t_n).$$
Taking the inverse limit of the maps $\kappa_n$,we obtain the bijection $\kappa_\infty: \cP_2 \to \mbZ_2$ of the path space $\cP_2$ onto the Cantor set $\mbZ_2$. Although $\mbZ_2$ is a group, the map $\kappa_\infty$ is only a homeomorphism, since $\cP_2$ does not have a group structure.

Recall that $N \cong \mbZ_2 \rtimes \mbZ_2^\times$, and $ \mbZ_2^\times$ acts on $\mbZ_2$ by \eqref{eq-normaction1}.
Recall \cite[Theorem 4.4.7]{RZ} that $\mbZ_2^\times \cong \mZ/2\mZ \times \mbZ_2$, where $\mZ/2\mZ$  is generated by $([-1]) \in \mbZ_2^\times$, where $[-1]$ denotes the equivalence class of $-1$ in $\mZ/2^n\mZ$ for $n \geq 1$, and the the second factor is generated by $([5]) \in \mbZ_2^\times$, where $[5]$ is the equivalence class of $5$ in $\mZ/2^n\mZ$ for $n \geq 1$. Denote these generators by $b$ and $c$ respectively. 

By \eqref{eq-normaction1}, $b$ acts on each $\mZ/2^n\mZ$ by taking an element of the group to its inverse, so $b$ has order $2$. Again by  \eqref{eq-normaction1}, $c$ acts on $\mZ/2^n\mZ$ by raising an element of the group to its $5$-th power. Since $b$ and $c$ are generators of a direct product of groups, they commute. This gives the relations in \eqref{eq-garith}. Using the identification of $V_n$ with $\mZ/2^n\mZ$, push forward the action of $b$ and $c$ to $\mZ/2^n\mZ$, then $b$ acts on $\mZ/2^n\mZ$ as multiplication by $-1$, and $c$ as multiplication by $5$. Note that both $b$ and $c$ fix the coset of $0$ in $\mZ/2^n\mZ$, therefore, $b$ and $c$ fix the word $x_n$ in $V_n$.

 For $n \geq 1$, denote by $G_n$ the isotropy group of the action of $G_{\rm arith}$ on $V_n$ at $x_n$. Since $G_{\rm geom} \subset G_{\rm arith}$ acts transitively on $V_n$, every coset in $G_{\rm arith}/G_n$ can be represented by a power of $a$, and so we can extend the maps \eqref{eq-kappa} to the bijections
  \begin{align}\label{eq-kapppa} \overline{\kappa}_n: G_{\rm arith}/G_n \to \mZ/2^n\mZ, \end{align}
for $n \geq 1$.  The coset of the identity in $G_{\rm arith}/G_n$ is mapped by \eqref{eq-kapppa} onto the coset of $0$ in $\mZ/2^n\mZ$. It follows that both $b$ and $c$ fix the coset of the identity in $G_{\rm arith}/G_n$, so we obtain that $\langle a^{2^n}, b, c\rangle \subset G_n$. 

Finally, suppose $g \in G_n$. Write $g$ as a word in $a,b,c$, that is, 
 $$g = a^{\alpha_1} c^{\gamma_1} b^{\beta_1} \cdots a^{\alpha_\ell} c^{\gamma_\ell} b^{\beta_\ell}, $$ 
where $\alpha_i, \beta_i,\gamma_i \in \mZ$ for $1 \leq i \leq \ell$. Using the relations in \eqref{eq-garith} we can rewrite $g$ in the form 
 $$g = a^{s} c^{k} b^{m}$$ 
for some $s, k, m \in \mZ$. Then $g \cdot x_n = x_n$ implies that $s$ is a power of $2^n$, and $g \in \langle a^{2^n},b,c \rangle$. We proved \eqref{eq-gngenerators-Q}. 
\endproof

We now compute the discriminant group $\cD_{\overline{x}}$ of the action of $G_{\rm arith}$ on $\cP_2$.
For $n \geq 1$, denote by $C_n$ the maximal normal subgroup of $G_n$ in $G$. 

\begin{lemma}\label{lemma-maxnormal}
Let $G_{\rm arith}$ be given by \eqref{eq-garith}, and $\{G_n\}_{n \geq 0}$ be a group chain given by \eqref{eq-gngenerators-Q}. Then for $n \geq 1$ the maximal normal subgroup of $G_n$ in $G_{\rm arith}$ is $C_n = \langle a^{2^n} , c^{2^n-2}\rangle \subset G_n$.
\end{lemma}

\proof The subgroup $C_n$ contains all elements of $G_n$ which fix every coset in $G_{\rm arith}/G_n$. So in particular $ \langle a^{2^n} \rangle \subset C_n$. 

Multiplication by $-1$ has order $2$ in $(\mZ/2^n\mZ)^\times$, and the generator $[5]$ has order $2^{n-2}$ in $(\mZ/2^n\mZ)^\times$ \cite[Theorem 5.44]{Rotman}, so the smallest power of $b$ which acts trivially on $G_{\rm arith}/G_n$ is $2$, and the smallest power of $c$ which acts trivially on $G_{\rm arith}/G_n$ is $2^{n-2}$. Thus $\langle a^{2^n}, b^2, c^{2^n-2} \rangle \subseteq C_n$.  

Applying the relations in \eqref{eq-garith}, every element of $G_{\rm arith}$ can be written down in the form $a^s b^k c^m$. We have to determine which compositions act trivially on $G_{\rm arith}/G_n$.

Since $b$ and $c$ fix the coset of the identity, for $s \ne 0 \mod 2^{n}$ the action of $a^s b^k c^m$ is non-trivial on the coset of the identity $eG_n$, so $a^s b^k c^m$ is not in $C_n$. 

Suppose $s = 0 \mod 2^n$, and $1 \leq m \leq 2^{n-2}-1$. We only need to consider the case when $k=1$, as $b^2=1$. Suppose $b c^m$ acts trivially on $G_{\rm arith}/G_n$, that is, for all $0 \leq t \leq 2^n-1$
  $$bc^m a^t G_n = a^t G_n .$$
Then $c^m a^t G_n = b a^t G_n$ for all $0 \leq t \leq 2^n-1$, and so $-1$ is in the subgroup of $(\mZ/2^n\mZ)^\times$ generated by $5$. But that is not the case  \cite[Theorem 5.44]{Rotman}, so $bc^m$ must act non-trivially on $G_{\rm arith}/G_n$ for all $1 \leq m \leq 2^{n-2}-1$. Thus we obtain that
   \begin{align}\label{eq-garithnc} C_n = \langle a^{2^n}, b^2,c^{2^{n-2}} \rangle = \langle a^{2^n}, c^{2^{n-2}} \rangle \subset G_n \subset G_{\rm arith}, \end{align}
 which finishes the proof of Lemma \ref{lemma-maxnormal}.
 \endproof  
   
 It follows from Lemma \ref{lemma-discr0} and Lemma  \ref{lemma-maxnormal} that for $n \geq 1$ there is an isomorphism
  \begin{align}\label{eq-quotientgroup-1}\lambda_n : G_n/C_n \to (\mZ/2^n\mZ)^\times \cong C_2 \times \mZ/2^{n-2}\mZ, \end{align}
where $C_2 = \{\pm 1\}$ is the multiplicative group of order $2$. The cosets in $G_n/C_n$ are represented by elements $b^k c^m$, where $k = 0,1$ and $0 \leq m \leq 2^{n-2}-1$, and 
  $$\lambda_n(b^k c^m C_n) = ((-1)^k, m + 2^{n-2}\mZ).$$
Taking the inverse limit we obtain that the discriminant group is a Cantor group, that is,
   $$\cD_{\overline{x}} = \lim_{\longleftarrow}\{G_{n+1}/C_{n+1} \to G_n/C_n\} \cong \mbZ_2^{\times}.$$  
 
We now compute the asymptotic discriminant. To do that, for each $m \geq 1$ we restrict the action to a clopen subset $U_m(x_m)$ of $\cP_2$ consisting of paths through the vertex $x_m = 0^m$ in $V_m$. The restricted action is that of the group $G_m$ given by \eqref{eq-gngenerators-Q}, and, associated to the action, there is a truncated group chain $\{G_n\}_{n \geq m}$. For each $m \geq 1$, we compute the discriminant group $\cD^m_{\overline{x}}$, and show that the natural map $\cD_{\overline{x}} \to \cD^m_{\overline{x}}$, described in more detail later, is an isomorphism.

As in Section \ref{subsec-groupchains} we have a homeomorphism
 $$U_m(x_m) \to {\ds G_{m,\infty}= \lim_{\longleftarrow}\{G_m/G_{n+1} \to G_m/G_{n}, \, n \geq m\}}.$$ 
Recall from the discussion before \eqref{eq-kapppa} that the cosets of $G_{\rm arith}/G_n$ are represented by the powers of $a$. Since $G_m = \langle a^{2^m}, b,c \rangle$, the cosets of $G_m/G_n$ are represented by powers of $a^{2^m}$, that is, the cosets are given by
   $$G_n, a^{2^m} G_n, a^{2 \cdot 2^m} G_n,\cdots, a^{(2^{n-m}-1)2^m}G_n.$$
Restricting the bijection \eqref{eq-kapppa} to the subgroup $G_m$, we obtain   
    \begin{align}\label{eq-kapppa-r} \overline{\kappa}_n: G_m/G_n \to 2^m\mZ/2^n \mZ \cong \mZ/ 2^{n-m}\mZ, \end{align}
 where $\overline{\kappa}_n(a^{s \cdot 2^m} G_n) = s \cdot 2^m + 2^n \mZ \mapsto s + 2^{n-m} \mZ$. As for  \eqref{eq-kapppa}, the action of $b$ and $c$ on $G_m/G_n$ pushes forward via \eqref{eq-kapppa-r}  to multiplication by $-1$ and by $5$ in $\mZ/2^{n-m}\mZ$ respectively. Therefore, the order of $c$ is equal to the order of multiplication by $5$ in $\mZ/2^{n-m} \mZ$, which is $2^{n-m-2}$. By an argument similar to the one in Lemma \ref{lemma-maxnormal} one obtains that the maximal normal subgroup of $G_n$ in $G_m$ is given by
   $$C^m_n = \langle a^{2^n},b^2,c^{2^{n-m-2}} \rangle = \langle a^{2^n},c^{2^{n-m-2}} \rangle \subset G_n \subset G_m,$$
and there is a group isomorphism
   \begin{align}\label{eq-quotientgroup-3}\lambda_{m,n}: G_n/C^m_n \to  (\mZ/2^{n-m}\mZ)^{\times} \cong C_2 \times \mZ/2^{n-m-2}\mZ .\end{align} 
The cosets in $G_n/C^m_n$ are represented by $b^k c^s$, where $k = 0,1$, and $0 \leq s \leq 2^{n-m-2}-1$, and 
  $$\lambda_n(b^k c^s C^m_n) = ((-1)^k, s + 2^{n-m-2}\mZ).$$
Then \eqref{eq-quotientgroup-3} implies that the discriminant group is a Cantor group, as we have
    $$\cD^m_{\overline{x}} = \lim_{\longleftarrow}\{G^{n+1}/C^m_{n+1} \to G_n/C^m_n,\, n \geq m\} \cong \mbZ_2^\times.$$
The last step is to construct the natural map from $\cD_{\overline{x}}$ to $\cD^m_{\overline{x}}$ and to show that it is an isomorphism. For that, consider the group homomorphisms 
  $$\psi_n:G_n/C_n \to G_n/C^{m}_{n},$$ 
given by coset inclusions. Combining these with the maps \eqref{eq-quotientgroup-1} and \eqref{eq-quotientgroup-3}, for $n \geq 1$ we obtain the following commutative diagram
\begin{align}\label{diag-1} \xymatrix{G_n/C_n \ar[rrrr]^{\psi_n:b^kc^sC_n \mapsto b^k c^sC^m_n} \ar[d]_{\lambda_n} & & & &G_n/C^m_n \ar[d]^{\lambda_{m,n}}  \\ C_2 \times \mZ / 2^{n-2} \mZ \ar[rrrr]_{\overline{\psi}_n:(\pm 1, s) \mapsto (\pm 1, s \, {\rm mod} \, 2^m)} &&& &C_2 \times \mZ/2^{n-m-2}\mZ}\end{align}
The diagram \eqref{diag-1} induces the maps of the inverse limits so that the following diagram is commutative
 \begin{align}\label{diag-2} \xymatrix{\cD_{\overline{x}} \ar[rrrr]^{\psi_\infty} \ar[d]_{\lambda_\infty} & & & &\cD^m_{\overline{x}} \ar[d]^{\lambda_{m,\infty}}  \\ \mbZ_2^{\times} \cong C_2 \times \mbZ_2 \ar[rrrr]_{\overline{\psi}_\infty} &&& &\mbZ_2^{\times} \cong C_2 \times\mbZ_2}\end{align}
 
where $\lambda_\infty$ and $\lambda_{m,\infty}$ are isomorphisms. So if $\overline{\psi}_\infty$ is injective, then $\psi_\infty: \cD_{\overline{x}} \to \cD^m_{\overline{x}}$ is injective.
 
So let $(s, (r_n)) \ne (t,(y_n)) \in C_2 \times \mbZ_2$. Since the maps $\overline{\psi}_n$ in \eqref{diag-1} are constant on the first component, if $s \ne t$, then $\overline{\psi}_\infty(s, (r_n)) \ne \overline{\psi}_\infty(t,(y_n))$ and we are done, so assume that $s = t$. Since the elements are distinct, then there exists $k \geq 1$ such that for all $n \geq k$ we have $r_n - y_n \ne 0 \mod 2^n$. Choose $n$ large enough so that $n-m> k$. Then
   $$r_n - y_n \ne 0 \mod 2^{n-m},$$
which shows that $\overline{\psi}_\infty(s, (r_n)) \ne \overline{\psi}_\infty(s,(y_n))$, and the map $\overline{\psi}_\infty$  is injective. This finishes the proof of Proposition \ref{prop-semidirect}.
\endproof

We now consider the case when the orbit of the critical point $c$ of $f(x)$ is strictly pre-periodic of length $2$. Recall from Theorem \ref{thm-1} that in this case the action of ${\rm Gal}_{\rm geom}(f)$ is stable with finite discriminant group. Recall that $\sigma$ denotes the non-trivial permutation of the set with two elements.

\begin{prop}\label{prop-semidirect-dihedral}
Let $K$ be a finite extension of $\mQ$, and let $f(x)$ be a quadratic polynomial with coefficients in the ring of integers of $K$. Suppose the orbit of the critical point $c$ of $f(x)$ is strictly pre-periodic with $\#P_c = 2$.Then the action of the profinite iterated arithmetic monodromy group ${\rm Gal}_{\rm arith}(f)$ on the space of paths $\cP_2$ of the tree $T$ of the solutions to $f^n(x) = t$, $n \geq 1$, is stable with infinite discriminant group.
\end{prop}

\proof Under the hypothesis of the proposition, ${\rm Gal}_{\rm geom}(f)$ is conjugate to the closure ${\rm CL}(\widetilde{H}_2)$ of the countable subgroup $\widetilde{H}_2$ generated by the elements $b_1 = \sigma$ and $b_2 = (b_1,b_2)$, see the proof of Theorem \ref{thm-1}. That is, there exists $w \in Aut(T)$ such that \cite[Theorem 3.4.1]{Pink2013}
 $${\rm CL}(\widetilde{H}_2) = w\, {\rm Gal}_{\rm geom}(f)\, w^{-1}.$$ 
Let $G_{\rm geom}$ be a dense subgroup which topologically generates ${\rm Gal}_{\rm geom}(f)$. 
Under the hypothesis of the proposition the element $w \in Aut(T)$ can be chosen in such a way that $w G_{\rm geom}w^{-1} = \widetilde{H}_2$ \cite[Section 3.6]{Pink2013}. So relabelling the vertices of the tree $T$ via the conjugating automorphism $w$, we may assume that  ${\rm Gal}_{\rm geom}(f) = {\rm CL}(\widetilde{H}_2) $ and $G_{\rm geom} = \widetilde{H}_2 $. 

The product $b_0 = b_1 b_2$ generates an odometer action, and by \cite[Proposition 3.1.9]{Pink2013} $G_{\rm geom} = \langle b_1 b_2 \rangle \rtimes \langle b_1 \rangle= \langle b_0 \rangle \rtimes \langle b_1 \rangle$ is infinite dihedral, with ${\rm Gal}_{\rm geom}(f) \cong \mbZ_2 \rtimes C_2 $, where $C_2 = \{ \pm 1\}$ is the multiplicative group of order $2$ generated by $b_1$, and $\mbZ_2$ is the group of dyadic integers, topologically generated by $b_0$.

Denote by $N$ the normalizer of ${\rm Gal}_{\rm geom}(f)$ in $Aut(T)$. By \cite[Proposition 3.5.2]{Pink2013} $N$  is isomorphic to the semi-direct product $\mbZ_2 \rtimes \mbZ_2^\times$ , where $\mbZ_2^\times \cong Aut(\mbZ_2)$ denotes the multiplicative group of the dyadic integers. This isomorphism is given explicitly by
   $$\mbZ_2 \rtimes \mbZ_2^\times \to N: (m, \ell) \mapsto b_0^{m} w_\ell,$$
 where $m = (m_i) \in \mbZ_2$, for $\ell = (\ell_i) \in \mbZ_2^\times$ we have 
   \begin{align*}z_\ell = (b_0^{\frac{1-\ell}{2}}z_\ell, z_\ell), \, w_\ell = b_0^{\frac{1 - \ell}{2}} z_\ell,\end{align*}
and for every $\ell \in \mbZ_2^\times$ we have
      \begin{align}\label{eq-normaction}w_\ell b_0 w_\ell^{-1} = b_0^\ell. \end{align}
In the rest of the proof we set $a = b_0$, so $a$ is an element generating the odometer action on $\cP_2$.
      
Since ${\rm Gal}_{\rm geom}(f)$ is a normal subgroup of ${\rm Gal}_{\rm arith}(f)$, we have ${\rm Gal}_{\rm arith}(f) \subseteq N \cong \mbZ_2 \rtimes \mbZ_2^\times$. By a similar argument to the one in Lemma \ref{lemma-discr0} we can identify the path space $\cP_2$ of the binary tree $T$ with the group of dyadic integers $\mbZ_2 $ generated by $a $. More precisely, since the group of the powers of $a$ acts transitively on the vertex sets $V_n$, for each $n \geq 1$, there is a bijection $\kappa_n: V_n \to \mZ/2^n\mZ$, given by \eqref{eq-kappa} and, taking the inverse limits as in Lemma \ref{lemma-discr0} we obtain a homeomorphism $\kappa_\infty: \cP_2 \to \mbZ_2$. 

Set $b = b_1$, then $b^2 = 1$, and we have 
  $$G_{\rm geom} = \langle a \rangle \rtimes \langle b \rangle = \langle a,b \mid b^2 = 1,\, bab^{-1} = a^{-1}\rangle.$$
The normalizer of the closure of the odometer action of the subgroup generated by $a$ in $Aut(T)$, which coincides with the normalizer of ${\rm Gal}_{\rm geom}(f)$ in $Aut(T)$, satisfies $N \cong \mbZ_2 \rtimes \mbZ_2^\times$. As in Lemma \ref{lemma-discr0}, we choose the equivalence classes $([-1])$ and $([5])$ as generators of $\mbZ_2^\times \cong C_2 \times \mbZ_2$, where $C_2 = \{ \pm 1\}$ is the multiplicative group of order $2$. Here $[-1]$ denotes the equivalence class of $-1$ in $\mZ/2^n\mZ$ for $n \geq 1$, and  $[5]$ denotes the equivalence class of $5$ in $\mZ/2^n\mZ$ for $n \geq 1$. Denote by $b$ and $c$ respectively the preimages of these generators in $N$. Then \eqref{eq-normaction} implies that $N$ is the closure of the action of the group
   \begin{align}\label{eq-H} H \cong \langle a,b,c \mid b^2 = 1, \, bab^{-1} = a^{-1}, \, cac^{-1} = a^5, \, bcb^{-1}c^{-1} = 1\rangle  \subset Aut(T), \end{align}
where the last relation follows from the fact that $b$ and $c$ correspond to the generators of the different factors of the product $\mbZ_2^\times  \cong C_2 \times \mbZ_2 $. 

Since $K$ is a finite extension of $\mQ$, by \cite[Corollary 3.10.6 (g)]{Pink2013} the quotient ${\rm Gal}_{\rm arith}(f)/{\rm Gal}_{\rm geom}(f)$ is a subgroup of finite index in the quotient $\mbZ_2 \rtimes \mbZ_2^\times / \mbZ_2 \rtimes \{ \pm 1\}$. This subgroup is generated by a power of $([5])$, and so the corresponding subgroup in $N/{\rm Gal}_{\rm geom}(f)$ is generated by a power $z = c^t$, for some $t \geq 1$. Since the order of $c$ is $2^{n-2}$, the order of $z = c^t$ is $2^{n-2}/t$ which implies that $t = 2^r$ for some $1\leq r \leq {n-2}$. Using the relation $ca = a^5c$, we then obtain that
  $$z a z^{-1} = c^t a c^{-t} = a^{5t} = a^{5\cdot {2^r}},$$
and so ${\rm Gal}_{\rm arith}(f)$ has a dense subgroup $G_{\rm arith}$ with a presentation
   \begin{align}\label{eq-garith-2} G_{\rm arith} \cong \langle a,b,c \mid b^2 = 1, \, bab^{-1} = a^{-1}, \,za z^{-1} = a^{5 \cdot {2^r}}, \, bzb^{-1}z^{-1} = 1\rangle  \subset Aut(T). \end{align}
As in Lemma \ref{lemma-discr0}, let $x_n = 0^n$ be the vertex labelled by a word of $0$'s in $V_n$, and let $\overline{x} = 0^\infty$ be the path containing the vertices $x_n$, $n \geq 1$. We choose $\overline{x}$ as our basepoint. Since the action of $c$ fixes $\overline{x}$, then the action of $z = c^{2^r}$ fixes $\overline{x}$. Arguing further similarly to Lemma \ref{lemma-discr0} we obtain that for each $n \geq 1$, the isotropy group of the action of $G_{\rm arith}$ at $x_n$ is given by
  \begin{align}\label{eq-gngenerators}G_n = \langle a^{2^n},b,z \rangle \subset G_{\rm arith}. \end{align}
Thus we obtain a group chain $\{G_n\}_{n \geq 1}$. 

Denote by $C_n$ the maximal normal subgroup of $G_n$ in $G_{\rm arith}$. Using that $z$ is a power of $c$, and arguing similarly to the proof of Lemma \ref{lemma-maxnormal} we obtain that 
  \begin{align}\label{eq-cngenerators}C_n = \langle a^{2^n},z^{2^{n-r-2}} \rangle \subset G_{\rm arith}. \end{align}
Then by an argument similar to the one after Lemma \ref{lemma-maxnormal} we obtain that for $n \geq 1$ we have an isomorphism
  $$ \lambda_n: G_n/C_n \to C_2 \times \mZ/2^{n-r-2}\mZ,$$
and, taking the inverse limits, we obtain an isomorphism $\lambda_\infty: \cD_{\overline{x}} \to \mbZ_2^\times$. Thus the discriminant group $\cD_{\overline{x}}$ is an infinite profinite group. 

We now compute the asymptotic discriminant. As in Proposition \ref{prop-semidirect}, for each $m \geq 1$ we restrict the action to a clopen subset $U_m(x_m)$ of $\cP_2$ consisting of paths through a vertex $x_m = 0^m$ in $V_m$. The restricted action is that of the group $G_m$ given by \eqref{eq-gngenerators} for $n = m$, and, associated to the action, there is a truncated group chain $\{G_n\}_{n \geq m}$. For each $m \geq 1$, we compute the discriminant group $\cD^m_{\overline{x}}$, and show that the natural map $\cD_{\overline{x}} \to \cD^m_{\overline{x}}$ is an isomorphism.

The argument is similar to the one in Proposition \ref{prop-semidirect}, with small adjustments for the fact that we are now looking at a specific finite index subgroup of $\mbZ_2 \rtimes \mbZ_2^\times$. 
As before, we have a homeomorphism
 $$U_m(x_m) \to {\ds G_{m,\infty}= \lim_{\longleftarrow}\{G_m/G_{n+1} \to G_m/G_{n}, \, n \geq m\}}.$$ 
 Since $\langle a \rangle \subset G_{\rm arith}$ acts transitively on $V_n$, every coset of $G_{\rm arith}/G_n$ is represented by a power of $a$, and so there is a bijection $\overline{\kappa}_n: G_{\rm arith}/G_n \to \mZ/2^n\mZ$ as in \eqref{eq-kapppa}. 
Since $G_m = \langle a^{2^m}, b,z \rangle$, the cosets of $G_m/G_n$ are represented by powers of $a^{2^m}$, and $G_m/G_n$ bijects under $\overline{\kappa}_n$ onto the subgroup $2^m\mZ/ 2^{n} \mZ \cong \mZ/2^{n-m} \mZ$.

By an argument similar to the one in Proposition \ref{prop-semidirect} for $c$, we obtain that the order of $z$ is equal to the order of multiplication by $5 \cdot {2^r}$ in $\mZ/2^{n-m} \mZ$, which is $2^{n-m -r - 2}$. Using an argument similar to the one in Lemma \ref{lemma-maxnormal} one obtains that the maximal normal subgroup of $G_n$ in $G_m$ is
   $$C^m_n = \langle a^{2^n},b^2,z^{2^{n-m-r-2}} \rangle = \langle a^{2^n},z^{2^{n-m-r-2}} \rangle \subset G_n \subset G_m,$$
and there is a group isomorphism
   \begin{align}\label{eq-quotientgroup-4}\lambda_{m,n}: G_n/C^m_n \to  (\mZ/2^{n-m-r}\mZ)^{\times} \cong C_2 \times \mZ/2^{n-m-r-2}\mZ .\end{align} 
The it follows by an argument similar to the one in  Proposition \ref{prop-semidirect}  that  
$$\cD^m_{\overline{x}} = \lim_{\longleftarrow}\{G^{n+1}/C^m_{n+1} \to G_n/C^m_n,\, n \geq m\} \cong \mbZ_2^\times, $$
and the natural maps $\psi_\infty: \cD_{\overline{x}} \to \cD_{\overline{x}}^m$, given by the inverse limits of coset inclusions $\psi_n:G_n/C_n \to G_n/C^{m}_{n}$ are isomorphisms. Thus the asymptotic discriminant of the action is stable with infinite discriminant group.
\endproof

This finishes the proof of Theorem \ref{thm-3}.


\section{Non-Hausdorff elements and the subgroups of profinite groups}\label{sec-subgroups}

In this section we prove Theorem \ref{thm-4}, that is, given an action of a profinite group $H_\infty$ on the path space of a spherically homogeneous tree $T$, we give a condition under which $H_\infty$ contains non-Hausdorff elements.

Let $T$ be a spherically homogeneous tree as in Section \ref{ex-treecylinder}, that is, the vertex set of $T$ is $\bigsqcup_{n \geq 0} V_n$, and there is a sequence of integers $(\ell_1,\ell_2,\cdots)$, called the spherical index of $T$, such that for $ n \geq 1$ every vertex in $V_{n-1}$ is connected  by edges to precisely $\ell_n$ vertices in $V_n$. We assume that $\ell_n \geq 2$ for $n \geq 1$. Recall that the set of paths $\cP$ in $T$ is a Cantor set by Section \ref{ex-treecylinder}. By assumption of Theorem \ref{thm-4} there is a profinite group ${\ds H_\infty = \lim_{\longleftarrow}\{H_{n+1} \to H_n\}}$ which acts on $\cP$ in such a way that for each $n \geq 1$ the restriction of the action of $H_\infty$ to $V_n$ is given by the action of $H_n$, and that $H_n$ acts transitively on $V_n$. By the assumption of Theorem \ref{thm-4} there is also a collection of finite groups $\{L_n\}_{n \geq 1}$ such that for each $n \geq 1$ the wreath product $L_n \rtimes L_{n-1} \rtimes \cdots \rtimes L_1 \subset H_n$. We are going to prove that $H_\infty$ contains a non-Hausdorff element. 

\emph{Proof of Theorem \ref{thm-4}}. If $g$ is non-Hausdorff, then there exists an infinite path $\overline{x}$, a descending collection of clopen neighborhoods $\{W_n\}_{n \geq 0}$ with $\bigcap_{n \geq 1} W_n = \{\overline{x}\}$ and, for each $n \geq 1$, a clopen subset $O_n \subset W_n$, such that $g(\overline{x}) = \overline{x}$, $g|O_n$ is the identity, while $g|W_n$ is not the identity homeomorphism. We now find such $\overline{x}$, $\{W_n\}_{n \geq 0}$ and $\{O_n\}_{n \geq 1}$ in $H_\infty$.

Let $x_n = 0^n$ denote a word of length $n$ consisting only of $0$'s, and let $\overline{x} = 0^\infty$ be an infinite sequence of $0$'s. Then $\overline{x}$ corresponds to a path in $\cP$ containing the vertices labelled by $x_n$, for $n \geq 1$. 

For $k \geq 1$, let $W_k = U_{2k}(x_{2k}) = U_{2k}(0^{2k})$, that is, $U_{2k}(x_{2k})$ contains all paths in $\cP$ through the vertex in $V_{2k}$ labelled by $0^{2k}$. Since $\ell_n \geq 2$ for all $n \geq 1$ by assumption, every $x_{2k} \in V_{2k}$ is connected to at least two vertices in $V_{2k+1}$, one labelled by $0^{2k}0$ and another labelled by $0^{2k}1$. We denote $O_k = U_{2k+1}(0^{2k}1)$. We will obtain a non-Hausdorff element $g$ by induction.

Denote by $R_n$ a set with $\ell_n$ elements, and note that there are bijections $b_n: R_1 \times \cdots \times R_n \to V_n$, see Section \ref{ex-treecylinder} for details. Now recall the definition of the wreath product. 
First, let $\cL_1 = L_1$, and  suppose that $\cL_n$ is defined. Let $L_{n+1}^{|V_n| }=\{f: V_n \to L_{n+1}\}$ be the set of all functions from $V_n$ to $L_{n+1}$. Then the group
 \begin{eqnarray} \label{wp-formula-1}\cL_{n+1} = L_{n+1}^{|V_{n}|}  \rtimes \cL_n := L_{n+1} \rtimes L_n \rtimes \cdots \rtimes L_{1}\end{eqnarray}
  acts on the product $V_n \times R_{n+1} $ by
    \begin{eqnarray}\label{wp-action-1} (f,s)(v_n,w) = (s(v_{n}), f(s(v_{n})) \cdot w). \end{eqnarray}
That is, for each $v_n \in V_n$ the element $(f,s)$ maps $ \{v_n\} \times R_{n+1}$ to $ \{s(v_n)\} \times R_{n+1}$, and then acts on $ \{s(v_n)\} \times R_{n+1}$ as $f(s(v_n)) \in L_{n+1}$.

For each $n \geq 1$, let $p_n$ be a non-trivial element of $L_n$. Denote by $1$ the identity in $L_n$.

Let $k=1$, and define $g$ to act trivially on the vertex sets $V_1$ and $V_{2}$. That is, $g|_{V_2} = g_2$, where $g_2 = (f_1,1) \in \cL_2$ and $f_1: V_1 \to L_2$ is the trivial constant function. 

Define $f_2:V_2 \to L_3$ by
  $$f_2(01) = p_3, \, f_2(w) = 1 \textrm{ for all }w \ne 01,$$
 and set $g_3 = g|_{V_3} = (f_2,g_2)$. Then $g_3 \in \cL_3$.
 
 For $k >1$, suppose $g_{2k-1} = g|_{V_{2k-1}} \in \cL_{2k-1}$ is defined. Let $f_{2k-1}: V_{2k-1} \to L_{2k}$ be the trivial function, that is, $f_{2k-1}(w) = 1$ for all $w \in V_{2k}$. Define $g_{2k} = (f_{2k-1},g_{2k-1})$, then $g_{2k} \in \cL_{2k}$.  Define 
   \begin{align*} f_{2k}(0^{2k-1}1) = p_{2k+1},  \, f_{2k}(w) = 1 \textrm{ for }w \ne 0^{2k-1}1, \, w \in V_{2k},\end{align*}
and set $g_{2k+1} = g|_{V_{2k+1}} = (f_{2k},g_{2k})$.   Then $g_{2k+1} \in \cL_{2k+1}$ by definition of the wreath product. Note that for all $k \geq 1$ the element $g$ is trivial on the sets $O_k = U_{2k+1}(0^{2k}1)$, and non-trivial on the sets $U_{2(k+1)}(0^{2k+1}1)$. Both $O_k$ and $U_{2(k+1)}(0^{2k+1}1)$ are subsets of $W_k = U_{2k}(0^{2k})$. 

By definition, $g$ acts non-trivially only on clopen sets of paths passing through vertices labelled by a word $0^{2k+1}1$ for some $k \geq 1$. Since the path $\overline{x} = 0^\infty$ clearly does not contain such a vertex, we have $g(\overline{x}) = \overline{x}$. Thus $g$, $\overline{x}$, $\{W_k\}_{k \geq 1}$ and $\{O_k\}_{k \geq 1}$ are as desired, and $g$ is a non-Hausdorff element in $H_\infty$. 
   
\endproof

\begin{ex}\label{ex-11}
{\rm 
The paper \cite{BFHJY} studied profinite iterated monodromy groups and arboreal representations for the polynomial
 $$f(z) = -2z^3 + 3z^2 $$
over number fields.  
They showed that for this polynomial, $G={\rm Gal}_{\rm geom}(f) = {\rm Gal}_{\rm arith}(f)$, and that  $G$ contains the infinite wreath product $[C_3]^\infty$, where $C_3$ is a cyclic group of order $3$. Theorem \ref{thm-4} implies that $G$ contains a non-Hausdorff element, and so the action of $G$ on the path space $\cP_3$ of the rooted tree $T$ has wild asymptotic discriminant.

The methods of \cite{BFHJY} were extended to a general class of normalized (dynamical) Belyi maps in \cite{BEK2018}. All such maps are post-critically finite rational maps, with orbits of all ramification points pre-periodic. The paper \cite{BEK2018} specified a class of dynamical Belyi maps for which ${\rm Gal}_{\rm geom}(f)$ is isomorphic to the infinite wreath product of the alternating group $A_d$ with itself, where $d$ is the degree of the map. Theorem \ref{thm-4} implies that in this case ${\rm Gal}_{\rm geom}(f)$ and ${\rm Gal}_{\rm arith}(f)$ contain non-Hausdorff elements, and so their actions have wild asymptotic discriminants.
}
\end{ex}



\begin{thebibliography}{10}

\bibitem{AHM2005}
{W.~Aitken, F.~Hajir and C.~Maire},
\newblock {\it Finitely ramified iterated extensions},
\newblock {\bf International Mathematics Research Notices}, 14:855--880, 2005.

\bibitem{ALC2009}
{J.~{\'A}lvarez L{\'o}pez and A.~Candel},
\newblock {\it Equicontinuous foliated spaces},
\newblock {\bf Math. Z.}, 263:725--774, 2009.
 
 \bibitem{ALM2016}
{J.~{\'A}lvarez L{\'o}pez and M.~Moreira Galicia},
\newblock {\it Topological {M}olino's theory},
\newblock {\bf Pacific. J. Math.}, 280:257--314, 2016.


\bibitem{Auslander1988}
{J.~Auslander},
\newblock {\bf Minimal flows and their extensions},
\newblock {North-Holland Mathematics Studies}, Vol. 153, {North-Holland Publishing Co., Amsterdam}, 1988.

\bibitem{BOERT1996}
{H.~Bass, M.~V.~Otero-Espinar, D.~Rockmore, C.~Tresser},
\newblock {\bf Cyclic Renormalization and Automorphism Groups of Rooted Trees},
\newblock LNM 1621, Springer 1996.

\bibitem{BN2008}
{L.~Bartholdi and V.~Nekrashevych}
\newblock{\it Iterated monodromy groups of quadratic polynomials, I},
\newblock{\bf Groups Geom. Dyn.}, 2:309-336, 2008.

\bibitem{BFHJY}
{R. L. Benedetto and X. Faber and B. Hutz and J. Juul and Y. Yasufuku},
\newblock{\it A large arboreal Galois representation for a cubic postcritically finite polynomial},
\newblock{\bf Res. Number Theory}, 3, Art. 29, 21pp., 2017.   
  
\bibitem{BEK2018}
{I.~I.~Bouw, \:O.~Ejder and V.~Karemaker},
\newblock{\it Dynamical Belyi maps and arboreal Galois groups}, arXiv:1811.10086.  

\bibitem{ClarkHurder2013}
{A.~Clark and S.~Hurder},
\newblock {\it Homogeneous matchbox manifolds},
\newblock {\bf Trans. Amer. Math. Soc.}, 365:3151--3191, 2013.

\bibitem{CHL2017}
{A.~Clark, S.~Hurder and O.~Lukina},
\newblock {\it Classifying matchbox manifolds},
\newblock {\bf Geometry and Topology}, 23(1):1--27, 2019.


\bibitem{DHL2016}
{J.~Dyer, S.~Hurder and O.~Lukina},
\newblock {\it The discriminant invariant of Cantor group actions},
\newblock {\bf Topology Appl.}, 208: 64--92, 2016.


\bibitem{DHL2017}
{J.~Dyer, S.~Hurder and O.~Lukina},
\newblock {\it Molino theory for matchbox manifolds},
\newblock {\bf Pacific J. Math.}, 289(1):91--151, 2017.


\bibitem{FO2002}
{R.Fokkink and L.Oversteegen},
\newblock {\it Homogeneous weak solenoids},
\newblock {\bf Trans. Am. Math. Soc.}, 354(9):3743--3755, 2002.



\bibitem{Ellis1960}
{R.~Ellis},
\newblock {\it A semigroup associated with a transformation group},
\newblock {\bf Trans. Amer. Math. Soc.}, 94:272--281, 1969.
     
\bibitem{Ellis1969}
{R.~Ellis},
\newblock {\bf Lectures on topological dynamics},
\newblock {W. A. Benjamin, Inc., New York}, 1969.

\bibitem{Ellis2014}
{D.~ Ellis and R.~Ellis},
\newblock {\bf Automorphisms and equivalence relations in topological dynamics},
\newblock {London Mathematical Society Lecture Note Series}, Vol. 412, Cambridge University Press, Cambridge, 2014.
    


\bibitem{EllisGottschalk1960}
{R.~Ellis and W.H.~Gottschalk},
\newblock {\it Homomorphisms of transformation groups},
\newblock {\bf Trans. Amer. Math. Soc.}, 94:258--271, 1969.


\bibitem{GZ2002}
{R.~I.~Grigorchuk and A.~ \.Zuk},
\newblock {\it On a torsion-free weakly branch group defined by a three state automaton},
\newblock {\bf International Journal of Algebra and Computations}, 12(1)-(2):223--246, 2002.


\bibitem{Haefliger1985}
{A.~Haefliger},
\newblock {\it Pseudogroups of local isometries}, in Differential Geometry (Santiago de Compostela, 1984), edited by L.A. Cordero,
\newblock {\bf Res. Notes in Math.}, 131:174--197, Boston, 1985.


\bibitem{HL2018}
{S.~Hurder and O.~Lukina},
\newblock {\it Orbit equivalence and classification of weak solenoids},
\newblock  arXiv: 1803.02098, to appear in {\bf Indiana Univ. Math. J.}

\bibitem{HL2017}
{S.~Hurder and O.~Lukina},
\newblock {\it Wild solenoids},
\newblock  {\bf Trans. Amer. Math. Soc.}, 371(7): 4493--4533, 2019.

\bibitem{Hughes2012}
{B.~Hughes},
\newblock {\it Trees, Ultrametrics, and Noncommutative Geometry},
\newblock {\bf Pure and Applied Mathematics Quarterly}, 8(1):221--312, 2012.

     
\bibitem{Jones2008}
 {R. Jones},
 \newblock{\it The density of prime divisors in the arithmetic dynamics of quadratic polynomials},
  {\bf J. London Math. Soc. (2)}, 78:523-544, 2008.       
     
\bibitem{Jones2013}
{R. Jones},
\newblock{\it Galois representations from pre-image trees: an arboreal survey}
\newblock in {\bf Actes de la {C}onf\'erence ``{T}h\'eorie des {N}ombres et
              {A}pplications''}, 107-136, 2013.
              
 \bibitem{Jones2015}
 {R. Jones},
 \newblock{\it Fixed-point-free elements of iterated monodromy groups},
  {\bf Trans. Amer. Math. Soc.}, 367(3):2023-2049, 2015.             
  
\bibitem{JKMT2015}
  {J.~Juul, P.~Kurlberg, K.~Madhu and T.~J.~Tucker},
     {\it Wreath products and proportions of periodic points},
   {\bf Int. Math. Res. Not. IMRN}, 13: 3944--3969, 2015.
 
\bibitem{Lubotzky1993}
{A.~Lubotzky},
\newblock {\it Torsion in profinite completions of torsion-free groups},
\newblock {\bf Quart. J. Math. Oxford Ser. (2)} 44:327--332, 1993.

\bibitem{Lukina2018}
{O.~ Lukina},
\newblock{Arboreal Cantor actions},
\newblock {\bf J. Lond. Math. Soc.}, 99(3): 678--706, 2019.

\bibitem{Nekr}
{V. Nekrashevych}, {\bf Self-similar groups}, Mathematical Survey and Monographs, 117, Americal Mathematical Society, Providence, RI, 2005.

\bibitem{Nekr2007}
{V. Nekrashevych}
\newblock{\it A minimal Cantor set in the space of $3$ generated groups},
\newblock{\bf Geom. Dedicata}, 124:153--190, 2007.


\bibitem{Nekr2018}
{V. Nekrashevych}
\newblock{\it Palindromic subshifts and simple periodic groups of intermediate growth},
\newblock{\bf Ann. of Math. (2)}, 187(3): 667--719, 2018.


\bibitem{Odoni1985}
{R. W. K. Odoni},
\newblock {\it The {G}alois theory of iterates and composites of polynomials},
\newblock {\bf Proc. London Math. Soc. (3)}, 51:385-414, 1985.

\bibitem{Pink2013}
{R.~Pink},
\newblock {\it Profinite iterated monodromy groups arising from quadratic polynomials},
\newblock {arXiv:1307.5678}.

\bibitem{Pink2013-2}
{R.~Pink},
\newblock {\it Profinite iterated monodromy groups arising quadratic morphisms with infinite post-critical orbits},
\newblock {arXiv:1309.5804}.


\bibitem{Renault1980}
{J.~Renault},
\newblock {\it A groupoid approach to {$C^*$}-algebras},
\newblock {\bf Lecture Notes in Math.},  vol. 793, 1980.

\bibitem{Renault2008}
{J.~Renault},
\newblock {\it Cartan subalgebras in {$C^*$}-algebras},
 \newblock {\bf Irish Math. Soc. Bull.}, 61:29--63, 2008.

  \bibitem{RZ} L. Ribes, and P. Zalesskii, \textbf{Profinite groups}, Springer-Verlag, Berlin 2000.


   \bibitem{Rotman} J.~Rothman, \textbf{An Introduction to the Theory of Groups}, GTM 148, 4th edition, Springer.


  \bibitem{Willard} S. Willard, \textbf{General Topology}, Dover Publications, Inc., Mineola, New York, 2004, xii+370pp.

  \bibitem{Wilson} J. S. Wilson, \textbf{Profinite groups}, London Mathematical Society Monographs, New Series, 19, 1998, xii+284pp.

\bibitem{Winkel1983}
{E.~Winkelnkemper},
\newblock {\it The graph of a foliation},
\newblock {\bf Ann. Global Ann. Geo.}, 1:51--75, 1983.

\end{thebibliography}
\end{document}